# Defender-Attacker-Target Game: Open-Loop Solution

Vladimir Turetsky* and Valery Y. Glizer†

Department of Applied Mathematics,
Ort Braude College of Engineering, Karmiel, Israel

**Abstract**

A defender-attacker-target problem with non-moving target is considered. This problem is modeled by a pursuit-evasion zero-sum differential game with linear dynamics and quadratic cost functional. In this game the pursuer is the defender, while the evader is the attacker. The objective of the pursuer is to minimize the cost functional, while the evader has two objectives: to maximize the cost functional and to keep a given terminal state inequality constraint. The open-loop saddle point solution of this game is obtained in the case where the transfer functions of the controllers for the defender and the attacker are of arbitrary orders. Then, this result is applied to the case of the first order controllers for the defender and the attacker. Numerical illustrating examples are presented.

---
*turetsky1@braude.ac.il
†valery48@braude.ac.il

# 1  Introduction

In the present paper, we study a defender-attacker-target problem. This problem is considered as a kind of interception problems, and it is modeled by a pursuit-evasion differential game.

The topic of interception problems studies an engagement between two moving vehicles. One of these vehicles tries to capture the other, while the other second vehicle tries to avoid such a capture. We call the first vehicle a pursuer, while the second is called an evader. Interception problems have been considered in the literature in various settings. If the evader's behaviour is predictable by the pursuer, this problem can be formulated in the framework of optimal control theory (see, e.g., (Bryson & Ho, 1975; Guelman & Shinar, 1984; Glizer, 1997) and references therein). If the evader's behavior is unpredictable for the pursuer, the interception problem can be formulated as a robust control problem (Glizer & Turetsky, 2012), or, in particular, as a finite horizon pursuit-evasion differential game (see e.g., (Bryson & Ho, 1975; Isaacs, 1965; Krasovskii & Subbotin, 1988; Shinar, Glizer, & Turetsky, 2013) and references therein). The objective of the pursuer is to minimize the miss distance (the closest separation between the vehicles), while the evader tries to maximize this miss distance and thus to avoid the capture.

In a defender-attacker-target problem, a defender can be considered as a pursuer, while an attacker can be considered as an evader. In such situations, the avoiding of the capture is not the main aim of the evader, because it tries not only to escape from the pursuer but also to hit an attacked object (target).

A defender-attacker-target problem was studied in a number of works in the literature. Thus in (Shaferman & Shima, 2010; Shima, 2011; Prokopov & Shima, 2013; Weiss, Shima, Castaneda, & Rusnak, 2017), the defender-attacker-target problem with a moving target was solved in the frame of control theory. In these papers it is assumed that either the attacker knows the controls of the defender and the target, or the defender and the target know the attacker's control. Based on these assumptions, the original problem was modeled by a linear-quadratic optimal control problem. Various types of a such control problem were proposed for this model: with and without hard controls' constraints, with and without terminal part of the cost functional, with and without integral part of the cost functional, with and without terminal state constraints. If the above mentioned assumptions are not made, the defender-attacker-target problem can be solved in the frame of differential games. Such an approach was proposed in the works (Lipman & Shinar, 1995; Perelman, Shima, & Rusnak, 2011; Rubinsky & Gutman, 2014; Liang, Peng, & Li, 2016; Casbeer, Garcia, & Pachter, 2017; Li & Cruz, 2011). In (Liang et al., 2016), the case of a moving target is considered. All the participants of the engagement have the simple kinematic equations in a plane. Two differential games are proposed to model the engagement from the defender and the target, and the attacker viewpoints. In each of these games, the heading angles of the participants are their controls. In the first game, the final range between the attacker and the target is the cost functional. The terminal equality constraint is imposed on the final range between the defender and the attacker. In the second game, the final range between the defender and attacker is the cost functional. The terminal equality constraint is imposed on the final range between the attacker and the target. Both games were solved analytically. In (Casbeer et al., 2017), two types of the defender-attacker-target problem, with one and two defenders, were analyzed. In each type of the problem, the participants of the engagement have the simple kinematic equations in a plane and their heading angles are their controls. The duration of the game is not prescribed, and the final time is the first time instant at which the separation of the defender (or at least one of the defenders) and the attacker equals zero. The terminal separation of the target and the attacker is the cost functional in the game. Detailed analysis of these games was carried out. The game of kind, modeling the original problem, also was analyzed in (Casbeer et al., 2017). In (Rubinsky & Gutman, 2014), the original defender-attacker-target problem was modeled by a differential game of kind, and this game was solved from the attacker viewpoint. In this game the behaviour of the engagement's participants is described by linear differential equations with zero-order transfer functions of the controllers, i.e., the lateral accelerations serve as the controls. These controls are subject to hard constraints. Two final time instants, for the defender-attacker engagement and attacker-target engagement, are given. At each of these time instants, a state inequality constraint is given. One of these constraints provides the avoidance of the capture of the attacker by the defender, while the other provides the capture of the target by the attacker. Conditions for the existence of the attacker's control, solving this game, were derived and the bang-bang attacker's state-feedback control itself was designed. In (Lipman & Shinar, 1995), the original defender-attacker-target problem was considered in the case of non-moving target. This problem was modeled by a differential game of degree. In this game the dynamics of the defender and the attacker are described by linear differential equations. The defender has the first-order controller, and its control is its lateral acceleration command. The attacker has the zero-order transfer functions of the controller, and its control is its lateral acceleration. Both controls are subject to hard constraints. The final time instants, for the defender-attacker engagement and attacker-target engagement, are given. The first of these time instance is the final time of the game. At this time instance, a state inequality constraint is given. This constraint is necessary for the capture of the target by the attacker. The cost functional in this game is the miss distance in the the defender-attacker engagement. The solution of this game yields bang-bang optimal state-feedback controls for the defender and the attacker. In (Li & Cruz, 2011), a classic finite horizon zero-sum linear-quadratic differential game (without terminal state constraints) was applied to model a defender-attacker-target problem in the following situations: (1) the target is non-moving; (2) the behaviour of the target is known to the attacker and the defender; (3) the target escapes from the attacker. In the general case of linear dynamics of the engagement's participants, the existence of the solution to the corresponding Riccati differential equation was assumed, while in the case of the simple motion dynamics of the participants, conditions for such an existence were derived. A linear-quadratic pursuit-evasion differential game (without terminal state constraints) also was used in (Perelman et al., 2011) as a model of a defender-attacker-target problem with a moving target. Each participant of the engagement has a controller with an arbitrary-order transfer function. For this game, the continuous-time analytic and discrete-time numeric state-feedback solutions were obtained. These solutions present the optimal evasion strategy of the target and the optimal pursuit strategy of the defender, as well as the optimal strategy of the attacker for pursuing the target and for evading the defender.

In the present paper, a defender-attacker-target problem with non-moving target is considered. In contrast with the works (Shaferman & Shima, 2010; Shima, 2011; Prokopov & Shima, 2013; Weiss et al., 2017), we assume that the attacker and the defender do not know the controls of each other. In such a case, it is impossible to use an optimal control problem as a model of a defender-attacker-target problem. Therefore, we model this problem by a differential game. Moreover, in contrast with (Rubinsky & Gutman, 2014), we consider the game of degree, and solve it from viewpoints of both, defender and attacker, players. In this game, the dynamics of the defender and the attacker are given by linear systems of differential equations. In contrast with the works (Liang et al., 2016; Casbeer et al., 2017; Rubinsky & Gutman, 2014; Lipman & Shinar, 1995), where the control of the attacker and/or the defender is either the heading angle, or the lateral acceleration, or the lateral acceleration command, we consider an essentially more reasonable from the practical viewpoint (and, therefore, much more complicated from the mathematical viewpoint) controllers of the players. Namely, the transfer functions of the defender's and attacker's controllers are of arbitrary and different, in general, orders with scalar controls. No hard constraints are imposed on these controls. The finite durations of the defender-attacker engagement and attacker-target engagement are given. The first of these durations is the game duration. The cost functional is a weighted sum of the square of the defender-attacker miss distance and the finite horizon integrals of the squares of the defender's and attacker's controls. The weights of the first two terms in the cost functional are positive, while the weight of the third term is negative. The cost functional is minimized by a proper choice of the defender's control and is maximized by a proper choice of the attacker's control. Thus, the considered game is a zero-sum linear-quadratic differential game. However, in contrast with the works (Li & Cruz, 2011; Perelman et al., 2011), where a classic zero-sum linear-quadratic differential game (without constraints) was used to model a defender-attacker-target problem, we impose a non-strict inequality constrain on the terminal state of the attacker. This inequality is necessary for the attacker to capture the target. Therefore, using this inequality in the game of the present paper, we convert this game to a considerably more adequate model of a defender-attacker-target problem.

It is important to note the following. If hard (geometric) constraints are imposed on the defender/attacker control, then the corresponding optimal control is as a rule of a bang-bang form (see e.g. (Lipman & Shinar, 1995; Rubinsky & Gutman, 2014)). Such a control is piecewise constant with the values of the upper and lower constraints' bounds. If a frequency of the switches between these values is high, then the control chattering occurs. The control chattering phenomenon is extremely undesirable in various applications. In the present paper, we replace the hard control constraints used in (Lipman & Shinar, 1995; Rubinsky & Gutman, 2014) with the soft control constraints, which are represented by integral quadratic terms of the defender and attacker controls in the cost functional. Such a replacement, along with properly chosen penalty coefficients for these terms, incentives the optimal controls to satisfy given hard constraints, while to be smooth and not to reach the boundaries of these hard constraints. The latter allows to avoid the highly undesirable control chattering. The comparison of a bang-bang control and a smooth control with respect to the chattering phenomenon can be found, e.g., in (Turetsky & Glizer, 2005, 2007).

Thus, in the present paper, the defender-attacker-target problem is modeled by the linear-quadratic differential game with the terminal state inequality constraint. We solve this game, using the penalty function approach.

The penalty function approach to solution of differential games was used in several works in the literature. Thus in (Heymann, Rajan, & Ardema, 1985), a smooth, monotonically decreasing and unbounded penalty function was used to relax state constraints (called the event constraints) in a zero-sum differential game. Then, based on this relaxation, a computational method was proposed for obtaining suboptimal state-feedback controls of the players. In (Choi & Tahk, 2000), a time-optimization zero-sum differential game with a terminal state inequality constraint was considered. The authors propose to included this constraint into the cost functional, using a discontinuous penalty function. The effectiveness of such an approach was tested numerically in the planar pursuit-evasion game where speed and heading angle of each player are subject to a control input. In (Carlson & Leitmann, 2012), a variational game with equality constraints was considered. This game is equivalent to an $N$-person nonzero-sum differential game with simple motion of the players and equality constraints imposed simultaneously on the state and the controls. These constraints are penalized yielding a family (with respect to a penalty parameter) of unconstrained games. Subject to proper assumptions, it was shown that there exists a sequence of the penalty parameter tending to infinity, along which open loop Nash equilibria of the unconstrained games converge to an open loop Nash equilibrium of the original game.

In the present paper, the penalty function method is applied to the solution of the original game in the following way. First, the decomposition of the game space into two non-intersecting regions is carried out. In the first region, the original game is equivalent to an unconstrained zero-sum linear-quadratic differential game. Conditions for the existence of its open-loop saddle point are derived and this saddle point is obtained. In the second region, subject to a proper general condition, the equivalence of the original game to a zero-sum linear-quadratic differential game with a terminal state equality constraint is shown. The latter game is solved by the penalty function method with a quadratic penalty function and a penalty parameter tending to infinity. Based on asymptotic analysis of the penalized constraint-free game and using proper assumptions, the existence of an open-loop saddle point solution of the equality constraint game in the second region is established, and such a solution is obtained. Then conditions in the terms of the original game's data, guaranteing that this solution also is a solution of the original game, are derived. Thus, the game of the present paper and the way to solve it by application of the penalty function method differ considerably from the games and their penalty function solutions in the above mentioned works.

The above made comparison of the present paper with the ones known in the literature clearly shows the significant novelty of the differential game considered in the present paper, as well as of the mathematical technique proposed for the game's solution and the obtained results.

The paper is organized as follows. In the next section, the pursuit-evasion differential game with the terminal state inequality constraint, modeled the defender-attacker-target problem, is rigorously formulated. Its reduction to a lower dimensional game is carried out. Stages of the reduced game's solution are shortly described. In Section 3, the state space of the reduced game is decomposed into two non-intersecting regions. It is shown that in the first region the reduced game is equivalent to the unconstraint game. The open-loop saddle point solution of this game is derived. In Section 4, subject to the general condition it is shown that the reduced game, considered in the second region, is equivalent to the game with the terminal state equality constraint. In Section 5, this equality constraint game is decomposed into two subgames, each of which is solved by application of the penalty function method yielding an open-loop saddle point solution. In Section 6, based on the solutions of the subgames, the open-loop saddle point solution of the equality constraint game in the second region is obtained. In Section 7, based on the results of Section 4, it is shown the concrete conditions (in the terms of the reduced game's data), guaranteeing that this solution also is a solution of the reduced game with the inequality constraint. Using the results of the previous sections, the case of the first order controllers for the defender and the attacker is treated in Section 8. Conclusions are placed in Section 9.

# 2 Problem Statement

## 2.1 Pursuit-evasion model

The engagement between the defender (*pursuer*) and the attacker (*evader*) is considered. The mathematical model of this scenario is based on the following assumptions: (i) the engagement takes place in a horizontal plane; (ii) both players have constant velocities; (iii) each player has a linear controller dynamics; (iv) the relative trajectory can be linearized with respect to the nominal collision geometry.

In Fig. 1, the schematic engagement geometry is depicted. The $X$-axis is the initial line of sight. The $Y$-axis is normal to the $X$-axis. The origin of the coordinate system is collocated with the target (T) position, which is also the initial position of the pursuer. The points $(x_p, y_p)$ and $(x_e, y_e)$ are current coordinates of the pursuer (P) and the evader (E), respectively; $V_p$, $V_e$ are their velocities; $a_p$, $a_e$ are their lateral accelerations; $\varphi_p, \varphi_e$ are the respective angles between the velocity vectors and the $X$-axis.

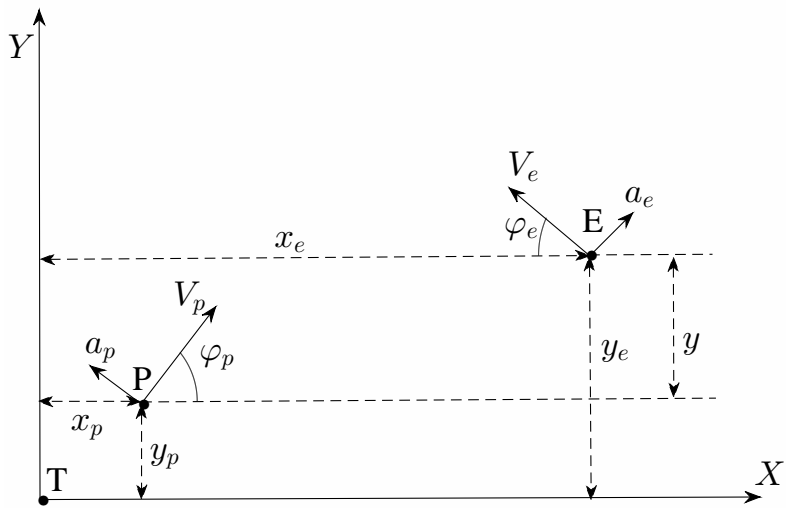

Figure 1: Interception geometry

Based on the small angles assumption (Glizer, Turetsky, Fridman, & Shinar, 2012; Shinar et al., 2013), the trajectories of the pursuer and the evader can be linearized with respect to the nominal collision geometry, leading to a constant closing velocity magnitude $V_c = |V_p| + |V_e|$. The final interception time $t_f$ can be easily calculated for any given initial range $r_0 = x_e(0)$: $t_f = r_0/V_c$. The value $x_f = |V_p|t_f = r_0 - |V_e|t_f$ is the distance in the $X$ direction between the evader and the target at the moment $t = t_f$.

Following the notation, introduced e.g. in (Shima, 2011), we assume that for $t \in [0, t_f]$, the controller dynamics of the pursuer and the evader are described by the equations

$$a_i = \bar{c}_i^T \bar{x}_i + \bar{d}_i u_i, \quad i = p, e, \tag{1}$$

$$\dot{\bar{x}}_i = \bar{A}_i \bar{x}_i + \bar{b}_i u_i, \quad \bar{x}_i(0) = [0]_{n_i \times 1}, \quad i = p, e, \tag{2}$$

where $\bar{x}_i$ is the state vector consisting of $n_i$ internal variables, $u_i$ is the scalar control, $i = p, e$, $[0]_{k \times m}$ denotes a zero $(k \times m)$-matrix. Note that $\bar{A}_i$ is an $(n_i \times n_i)$-matrix, $\bar{b}_i, \bar{c}_i \in R^{n_i}$, $\bar{d}_i \in R$ $(i = p, e)$.

For example, if the pursuer and/or the evader has the first-order strictly proper dynamics, then $n_i = 1$, $\bar{A}_i = -1/\tau_i$, $\bar{b}_i = 1/\tau_i$, $\bar{c}_i = 1$, $\bar{d}_i = 0$, where $\tau_i$ is the player's controller time constant; the control $u_i$ is the lateral acceleration command.

Due to (1) – (2), the motion of the pursuer and the evader are described by the differential equations

$$\dot{X}_i = A_i X_i + B_i u_i, \quad i = p, e \tag{3}$$

where the state vector is of the column block form:

$$X_i = [y_i, \dot{y}_i, \bar{x}_i^T]^T \triangleq [X_{i_1}, X_{i_2}, X_{i_3}, ..., X_{i_{n_i+2}}]^T, \quad i = p, e, \tag{4}$$

$$A_i = \begin{bmatrix} 0 & 1 & [0]_{1 \times n_i} \\ 0 & 0 & \bar{c}_i^T \\ [0]_{n_i \times 1} & [0]_{n_i \times 1} & \bar{A}_i \end{bmatrix}, \quad B_i = \begin{bmatrix} 0 \\ \bar{d}_i \\ \bar{b}_i \end{bmatrix}, \quad i = p, e. \tag{5}$$

Let us define the column block-form state vector

$$x = [X_p^T, X_e^T]^T = [y_p, \dot{y}_p, \bar{x}_p^T, y_e, \dot{y}_e, \bar{x}_e^T]^T. \tag{6}$$

Then, due to (3) – (6) and the small angles assumption, the system dynamics for $t \in [0, t_f]$ is described by the linear differential equation

$$\dot{x} = Ax + Bu_p + Cu_e, \tag{7}$$

where

$$A \triangleq \begin{bmatrix} A_p & [0]_{(n_p+2)\times(n_e+2)} \\ [0]_{(n_e+2)\times(n_p+2)} & A_e \end{bmatrix}, \tag{8}$$

$$B \triangleq \begin{bmatrix} B_p \\ [0]_{(n_e+2)\times 1} \end{bmatrix}, \quad C \triangleq \begin{bmatrix} [0]_{(n_p+2)\times 1} \\ B_e \end{bmatrix}. \tag{9}$$

The initial condition is $x(0) = [0, |V_p|\varphi_p^0, [0]_{1\times n_p}, 0, |V_e|\varphi_e^0, [0]_{1\times n_e}]^T$, where $\varphi_e^0$ and $\varphi_p^0$ are the (small) initial values of the aspect angles $\varphi_e$ and $\varphi_p$.

The objective of the pursuer is to minimize the cost functional

$$J(u_p(\cdot), u_e(\cdot)) = (y_e(t_f) - y_p(t_f))^2 + \alpha \int_0^{t_f} u_p^2(t)dt - \beta \int_0^{t_f} u_e^2(t)dt, \tag{10}$$

where $\alpha, \beta > 0$ are penalties for the players' controls.

**Remark 1** *The terminal part of the cost functional (10) is the squared miss distance, while the integral part represents the soft control constraints. In contrast with (Lipman & Shinar, 1995; Rubinsky & Gutman, 2014), we do not impose the geometric (hard) constraints on the controls. However, a proper choice of the coefficients $\alpha$ and $\beta$ can provide the validity of feasible geometric constraints for the players' optimal controls. Along with this, the structure of these controls differs considerably from the structure of the optimal controls in (Lipman & Shinar, 1995; Rubinsky & Gutman, 2014). Namely, the optimal controls of (Lipman & Shinar, 1995; Rubinsky & Gutman, 2014) are piecewise constant with the values of the upper and lower constraints' bounds, while the optimal controls obtained in the present paper are smooth and do not reach the boundaries of the feasible geometric constraints.*

The evader has a twofold objective. Its first objective is to maximize (10). The second objective is to be capable to reach the target at $t = t_f + t_c$, where $t_c$ is such that $x_e(t_f + t_c) = 0$: $t_c = x_f/|V_e| = \nu t_f$, $\nu \triangleq |V_p|/|V_e|$. In order to formulate the evader's second objective, let us calculate the $Y$-separation between the evader and the target for $t = t_f + t_c$. Due to (3) for $i = e$,

$$y_e(t_f + t_c) = X_{e1}(t_f + t_c) =$$

$$D_e \left( \Phi_e(t_f + t_c, t_f) X_e(t_f) + \int_{t_f}^{t_f+t_c} \Phi_e(t_f + t_c, t) B_e u_e(t) dt \right), \tag{11}$$

where
$$D_e = \left[1, [0]_{1\times(n_e+1)}\right], \tag{12}$$

$$\Phi_e(t_f + t_c, t) = \left\{\phi_{e_{kj}}(t_f + t_c, t)\right\}\Big|_{k,j=1}^{n_e+2}, \tag{13}$$

is the transition matrix of the homogeneous system, corresponding to (3) for $i = e$. This matrix satisfies

$$\frac{d}{dt}\Phi_e(t_f + t_c, t) = -\Phi_e(t_f + t_c, t)A_e, \ \Phi_e(t_f + t_c, t_f + t_c) = I_{n_e+2}. \tag{14}$$

For $t \in [t_f, t_f + t_c]$, we make a practically justified assumption that the evader's control is bounded:

$$|u_e(t)| \leq a_e^{\max}, \ a_e^{\max} > 0, \tag{15}$$

where $-a_e^{\max}$ and $a_e^{\max}$ are the lower and the upper bounds of the evader's control, respectively. Since $y_e(t_f + t_c) = 0$, by virtue of (11) and (15),

$$|D_e\Phi_e(t_f + t_c, t_f)X_e(t_f)| \leq \mu_e a_e^{\max}, \ \mu_e \triangleq \int_{t_f}^{t_f+t_c} |D_e\Phi_e(t_f + t_c, t)B_e|\,dt. \tag{16}$$

The pursuit-evasion differential game for the system (7) with the cost functional (10) and the evader's terminal constraint (16) is called the Original Game (OG).

**Remark 2** *A linear-quadratic differential game without hard control constrains is an acceptable model for various interception problems (see, e.g., (Bryson & Ho, 1975; Ben-Asher & Yaesh, 1998)). In this paper, the interception problem is mainly considered from the pursuer's viewpoint. Therefore, the absence of the evader's geometric control constraint for $t \in [0, t_f]$ in the OG represents the worst case from the pursuer's point of view, i.e., a hard geometric constraint is replaced by a soft constraint (the integral term of the cost functional (10)). From the other hand, the absence of the geometric constraint for the pursuer allows designing a more implementable pursuer's control, namely, a smooth open-loop control. Moreover, by a proper choice of the penalty coefficients $\alpha$ and $\beta$, we can provide a prescribed geometric constraint of the optimal $u_p(t)$. It should be also noted that in a real-life situation, the pursuer does not know the evader's behavior, and, consequently, cannot guarantee the capture (i.e., the condition $y_e(t_f) = y_p(t_f)$). Therefore, the pursuer is forced to play a game as the most effective tool for designing the control.*

## 2.2 Game reduction

Let us formulate the system, describing the relative motion between the evader and the pursuer in the direction normal to the initial line-of-sight (the $Y$-axis direction). Its state vector

$$X_{ep} = [y_e - y_p, \dot{y}_e - \dot{y}_p, x_p^T, x_e^T]^T \triangleq$$

$$[X_{ep_1}, X_{ep_2}, X_{ep_3}, ..., X_{ep_{n_p+2}}, X_{ep_{n_p+3}}, ..., X_{ep_{n_p+n_e+2}}]^T \quad (17)$$

satisfies the differential equation

$$\dot{X}_{ep} = A_{ep} X_{ep} + B_{ep} u_p + C_{ep} u_e, \quad t \in [0, t_f], \quad (18)$$

where $A_{ep} = \begin{bmatrix} 0 & 1 & [0]_{1 \times n_p} & [0]_{1 \times n_e} \\ 0 & 0 & -c_p^T & c_e^T \\ [0]_{n_p \times 1} & [0]_{n_p \times 1} & A_p & [0]_{n_p \times n_e} \\ [0]_{n_e \times 1} & [0]_{n_e \times 1} & [0]_{n_e \times n_p} & A_e \end{bmatrix}$,

$$B_{ep} = \begin{bmatrix} 0 \\ -\bar{d}_p \\ \bar{b}_p \\ [0]_{n_e \times 1} \end{bmatrix}, \quad C_{ep} = \begin{bmatrix} 0 \\ \bar{d}_e \\ [0]_{n_p \times 1} \\ \bar{b}_e \end{bmatrix}. \quad (19)$$

Let us introduce two new scalar state variables. The first is

$$z(t) \triangleq D_{ep} \Phi_{ep}(t_f, t) X_{ep}(t), \quad t \in [0, t_f], \quad (20)$$

where

$$D_{ep} = \begin{bmatrix} 1, [0]_{1 \times (n_p + n_e + 1)} \end{bmatrix}, \quad (21)$$

$\Phi_{ep}(t_f, t)$ is the transition matrix of the homogeneous system, corresponding to (18). The variable $z(t)$ is the zero-effort miss distance (ZEM) in the engagement of the evader and the pursuer: if at some moment $t = t_1$ both players stop to control ($u_p = u_e = 0$ for $t \in [t_1, t_f]$), then the value $y_e(t_f) - y_p(t_f)$ (miss distance) is equal to $z(t_1)$. Due to (17) and (20),

$$z(t_f) = y_e(t_f) - y_p(t_f). \quad (22)$$

The second scalar state variable is

$$w(t) \triangleq D_e \Phi_e(t_f + t_c, t) X_e(t), \quad (23)$$

where $D_e$ and $\Phi_e$ are given by (12) – (14). The variable $w(t)$ is the ZEM in the engagement between the evader and the static target.

By using (3), (14), (18) and the properties of the transition matrix, one obtains the differential equations for the new state variables $z$ and $w$ in the interval $[0, t_f]$:

$$\dot{z} = h_p(t)u_p + h_e(t)u_e, \quad z(0) = z_0, \tag{24}$$

$$\dot{w} = g_e(t)u_e, \quad w(0) = w_0, \tag{25}$$

where $g_e(t) = D_e\Phi_e(t_f+t_c, t)B_e$, $z_0 = t_f(|V_e|\varphi_e^0 - |V_p|\varphi_p^0)$, $w_0 = (t_f+t_c)|V_e|\varphi_e^0$,

$$h_p(t) = D_{ep}\Phi_{ep}(t_f, t)B_{ep}, \quad h_e(t) = D_{ep}\Phi_{ep}(t_f, t)C_{ep}. \tag{26}$$

Due to (22), the cost functional (10) can be rewritten as

$$J(u_p(\cdot), u_e(\cdot)) = z^2(t_f) + \alpha \int_0^{t_f} u_p^2(t)dt - \beta \int_0^{t_f} u_e^2(t)dt. \tag{27}$$

Due to (23), the constraint (16) becomes

$$|w(t_f)| \leq \mu_e a_e^{\max}. \tag{28}$$

In what follows, we deal with the pursuit-evasion differential game for the system (24) – (25) with the cost functional (27) and the terminal evader's constrain (28). We call this game the Reduced Game (RG), because its state vector $(z, w)^T$ is only two-dimensional, whereas the state vector $x$ of the original game has the dimension $n_p + n_e + 4$. In this game, the set of admissible pursuer's controls is

$$U_p = \{u_p(\cdot) \in L_2[0, t_f]\}, \tag{29}$$

where $L_2[0, t_f]$ is the space of square integrable functions defined on $[0, t_f]$, whereas the set of admissible evader's controls is

$$U_e = \left\{ u_e(\cdot) \in L_2[0, t_f] : \left| w_0 + \int_0^{t_f} g_e(t)u_e(t)dt \right| \leq \mu_e a_e^{\max} \right\}. \tag{30}$$

## 2.3 Solution stages

The RG solution is carried out in several stages. First, we consider the corresponding game without the constraint (28) (the Unconstrained Reduced Game, URG). We obtain the condition, guaranteeing that the URG open-loop saddle point solution satisfies the constraint (28). In this case, the solution of the URG is also the solution of the RG. Then, we show that if the above mentioned condition is not valid, the RG open-loop saddle point solution yields the equality terminal constraint

$$|w(t_f)| = \mu_e a_e^{\max}. \tag{31}$$

The game with the constraint (31) is called the Equality constraint Reduced Game (ERG). For this game we establish a general sufficient condition, subject to which its open-loop saddle point solution becomes a solution in the RG. Then, the open-loop saddle point solutions of two ERG subgames: ERG$^+$ and ERG$^-$ for the constraints $w(t_f) = \mu_e a_e^{\max}$ and $w(t_f) = -\mu_e a_e^{\max}$, respectively, are derived. Based on the solutions of the ERG$^+$ and the ERG$^-$, conditions for the existence of the ERG open-loop saddle point solution are established, and the saddle point itself is derived. Finally, we show that the saddle point of the ERG becomes the saddle point of the RG without any additional condition.

## 3 Unconstrained Game Solution

In the URG, the differential equation (25) can be dropped from the system dynamics. Once having the optimal control $u_e^0(t)$, the value of $w(t_f)$ is obtained straightforwardly by solution of (25) for $u_e = u_e^0(t)$. Knowing $w(t_f)$, we can derive the conditions guaranteeing that (28) is satisfied.

The admissible controls in this game for both players are square-integrable functions in the interval $[0, t_f]$. The solution of a linear-quadratic differential game is well known (see, e.g., (Bryson & Ho, 1975; Bernhard, 2015), and references therein).

Denote

$$s \triangleq 1 + \frac{1}{\alpha} \int_0^{t_f} h_p^2(t) dt - \frac{1}{\beta} \int_0^{t_f} h_e^2(t) dt. \tag{32}$$

**Proposition 1** *If*

$$\beta > \int_0^{t_f} h_e^2(t)dt, \qquad (33)$$

*then $s \neq 0$ and the pair*

$$u_p^0(t) = -\frac{h_p(t)z_0}{\alpha s}, \quad u_e^0(t) = \frac{h_e(t)z_0}{\beta s} \qquad (34)$$

*constitutes the URG saddle point.*

**Proof.** The expressions (34) directly follow from (Bernhard, 2015, Eq. (3)). Due to (Bernhard, 2015, Theorem 1), the pair (34) is the URG saddle point if the solution of the differential Riccati equation

$$\dot{P} = -\frac{1}{\beta}h_e^2(t)P^2, \quad P(t_f) = 1, \qquad (35)$$

has no conjugate points in the interval $[0, t_f]$. The solution of (35) is $P(t) = \left[1 - \frac{1}{\beta}\int_t^{t_f} h_e^2(\xi)d\xi\right]^{-1}$. Thus, this solution has no conjugate points if the condition (33) holds. □

**Remark 3** *Subject to the condition (33), $s > 0$, and, consequently, the open-loop saddle-point controls (34) are feasible for $t \in [0, t_f]$.*

Due to (25) and (34), $w(t_f) = w_0 + \int_0^{t_f} g_e(t)u_e^0(t) = w_0 + az_0$, where

$$a \triangleq \frac{1}{\beta s}\int_0^{t_f} h_e(t)g_e(t)dt. \qquad (36)$$

Thus, the condition (28) can be rewritten as $|w_0 + az_0| \leq \mu_e a_e^{\max}$. For fixed $\nu$, $a_e^{\max}$, $\alpha$ and $\beta$, the strict version of this defines the set $\Omega$ in the plane $(z_0, w_0)$. It is an infinite strip between the straight lines $w_0 = -az_0 \pm \mu_e a_e^{\max}$:

$$\Omega = \left\{(z_0, w_0) : |w_0 + az_0| < \mu_e a_e^{\max}\right\}. \qquad (37)$$

# 4 Reduced Game for $(z_0, w_0) \notin \Omega$

**Lemma 1** *Let $(u_p^*(\cdot), u_e^*(\cdot))$ be a saddle point in the RG with $(z_0, w_0) \notin \Omega$. Then, the solution $w(t)$ generated by $u_e^*(\cdot)$, satisfies the equality (31).*

**Proof.**

Since $(u_p^*(\cdot), u_e^*(\cdot))$ is the RG saddle point, then the saddle point inequality is valid:

$$J(u_p^*(\cdot), u_e(\cdot)) \leq J(u_p^*(\cdot), u_e^*(\cdot)) \leq J(u_p(\cdot), u_e^*(\cdot)), \ \forall \ u_p(\cdot) \in U_p, u_e(\cdot) \in U_e, \tag{38}$$

where the sets of the admissible pursuer's and evader's controls are given by (29) and (30), respectively.

The left-hand inequality in (38) means that $u_e^*(t)$ is the solution of the optimal control problem with the following data: (i) the dynamics (24) – (25) with $u_p(t) = u_p^*(t)$, (ii) the cost functional (27) with $u_p(t) = u_p^*(t)$ to be maximized by $u_e(t)$, and (iii) the terminal constraint (28). Then, due to (Vasilyev, 1988), there exist the functions $\lambda_z(t)$, $\lambda_w(t)$ and the numbers $\lambda_1^f \geq 0$, $\lambda_2^f \geq 0$, such that

$$\dot{\lambda}_z = 0, \quad \lambda_z(t_f) = 2z(t_f), \tag{39}$$

$$\dot{\lambda}_w = 0, \quad \lambda_w(t_f) = \lambda_1^f - \lambda_2^f, \tag{40}$$

$$\lambda_1^f(w(t_f) - \mu_e a_e^{\max}) = 0, \quad \lambda_2^f(w(t_f) + \mu_e a_e^{\max}) = 0, \tag{41}$$

and

$$u_e^*(t) = \frac{1}{2\beta}\left[h_e(t)\lambda_z(t) + g_e(t)\lambda_w(t)\right], \tag{42}$$

where $z(t)$ and $w(t)$ are the solutions of (24) – (25), generated by the pair $(u_p^*(\cdot), u_e^*(\cdot))$. Note that (39) yields

$$\lambda_z(t) \equiv 2z(t_f). \tag{43}$$

Now, let us assume that the statement of the lemma is wrong, i.e the solution $w(t)$ generated by $u_e^*(\cdot)$, satisfies

$$w(t_f) < \mu_e a_e^{\max}, \quad w(t_f) > -\mu_e a_e^{\max}. \tag{44}$$

In this case, the complementary slackness condition (41) leads to $\lambda_1^f = \lambda_2^f = 0$, which, due to (40), yields

$$\lambda_w(t) \equiv 0. \tag{45}$$

Due to (43) and (45),
$$u_e^*(t) = \frac{h_e(t)z(t_f)}{\beta}. \tag{46}$$

Now, let us specify the strategy $u_p^*(t)$. The right-hand side inequality in (38) means that $u_p^*(t)$ is the solution of the optimal control problem with the following data: (i) the dynamics (24) with $u_e(t) = u_e^*(t)$, and (ii) the cost functional (27) with $u_e(t) = u_e^*(t)$ to be minimized by $u_p(t)$. The solution of this standard linear-quadratic control problem without terminal constrain is (Vasilyev, 1988)
$$u_p^*(t) = -\frac{h_p(t)z(t_f)}{\alpha}. \tag{47}$$

Note that $z(t_f)$ in (46) and (47) is the same: it is the terminal value of the solution of (24), generated by the pair $(u_p^*(\cdot), u_e^*(\cdot))$. By substituting (46) and (47) into (24) and solving the resulting equation: $z(t_f) = z_0 - \frac{z(t_f)}{\alpha}\int_0^{t_f} h_p^2(t)dt + \frac{z(t_f)}{\beta}\int_0^{t_f} h_e^2(t)dt$, yielding

$$z(t_f) = \frac{z_0}{s}, \tag{48}$$

where the function $s$ is given by (32). By substituting (48) into (46) and (47), one can see that, subject to the assumption (44), the RG saddle point $(u_p^*(\cdot), u_e^*(\cdot))$ coincides with the URG saddle point (34). Therefore, by definition (37) of the set $\Omega$, this means that $(z_0, w_0) \in \Omega$ which contradicts with the assumption of the lemma. Thus the assumption (44) is wrong. This proves the lemma. □

Due to Lemma 1, for the case $(z_0, w_0) \notin \Omega$, an open-loop saddle point solution of the RG is an open-loop saddle point solution of the ERG. The latter consists of the dynamics equations (24) – (25), the cost functional (27) and the terminal evader's equality constraint (31). In this game, the set of admissible pursuer's controls is $U_p = \{u_p(\cdot) \in L_2[0, t_f]\}$, whereas the set of admissible evader's controls is

$$\bar{U}_e = \left\{ u_e(\cdot) \in L_2[0, t_f] : \left| w_0 + \int_0^{t_f} g_e(t)u_e(t)dt \right| = \mu_e a_e^{\max} \right\}. \tag{49}$$

Now, let us derive a condition, subject to which an open-loop saddle point solution $(\bar{u}_p^*(\cdot), \bar{u}_e^*(\cdot))$ of the ERG is an open-loop saddle point solution of the

RG. For this purpose, we consider the optimal control problem, consisting of the equations of dynamics

$$\dot{z} = h_p(t)\bar{u}_p^*(t) + h_e(t)\bar{u}_e, \quad z(0) = z_0, \tag{50}$$

$$\dot{w} = g_e(t)\bar{u}_e, \quad w(0) = w_0, \tag{51}$$

the state constraint (28) and the performance index

$$\bar{J}_e(\bar{u}_p^*(\cdot), \bar{u}_e(\cdot)) \triangleq z^2(t_f) - \beta \int_0^{t_f} \bar{u}_e^2(t)dt \to \max_{\bar{u}_e(\cdot) \in U_e}. \tag{52}$$

**Lemma 2** *Let $(\bar{u}_p^*(\cdot), \bar{u}_e^*(\cdot))$ be a saddle point in the ERG with $(z_0, w_0) \notin \Omega$. Then, the pair $(\bar{u}_p^*(\cdot), \bar{u}_e^*(\cdot))$ also is a saddle point in the RG with the same initial position $(z_0, w_0)$ if and only if $\bar{u}_e(t) = \bar{u}_e^*(t)$ is an optimal control in the problem (50) – (52),(28).*

**Proof.** *Necessity.* Let $(\bar{u}_p^*(\cdot), \bar{u}_e^*(\cdot))$ be a saddle point in the RG. Then, the inequality

$$J(\bar{u}_p^*(\cdot), u_e(\cdot)) \leq J(\bar{u}_p^*(\cdot), \bar{u}_e^*(\cdot)) \leq J(u_p(\cdot), \bar{u}_e^*(\cdot)), \ \forall \ u_p(\cdot) \in U_p, u_e(\cdot) \in U_e. \tag{53}$$

is satisfied. The left-hand side of this inequality directly implies that $\bar{u}_e(t) = \bar{u}_e^*(t)$ is an optimal control in the problem (50) – (52), (28). Thus, the necessity is proven.

*Sufficiency.* Since $(\bar{u}_p^*(\cdot), \bar{u}_e^*(\cdot))$ is a saddle point in the ERG, then it satisfies the inequality

$$J(\bar{u}_p^*(\cdot), u_e(\cdot)) \leq J(\bar{u}_p^*(\cdot), \bar{u}_e^*(\cdot)) \leq J(u_p(\cdot), \bar{u}_e^*(\cdot)), \ \forall \ u_p(\cdot) \in U_p, u_e(\cdot) \in \bar{U}_e. \tag{54}$$

To prove the sufficiency, we should show the fulfilment of the inequality (53). The fulfilment of the right-hand inequality in (53) directly follows from the right-hand side inequality in (54). The fulfilment of the left-hand inequality in (53) immediately follows from the assumption that $\bar{u}_e^*(t)$ is an optimal control in the problem (50) – (52), (28). Thus, the sufficiency is proven, which completes the proof of the lemma. □

**Remark 4** *Lemmas 1 and 2 imply the following. If the assumption of Lemma 2 with respect to the optimal control problem (50) – (52), (28) is satisfied, then the RG and the ERG with the same initial position $(z_0, w_0) \notin \Omega$* *are equivalent to each other in the sense of the open-loop saddle point solutions. In the subsequent sections, we obtain an open-loop saddle point in the ERG for $(z_0, w_0) \notin \Omega$. Then, we show that this saddle point satisfies the assumption of Lemma 2 with respect to the optimal control problem (50) – (52), (28), i.e., the obtained ERG saddle point becomes a saddle point in the RG.*

## 5 Saddle Point Solutions of Equality Constraint Subgames for $(z_0, w_0) \notin \Omega$

For the solution of these subgames $\text{ERG}^+$ and $\text{ERG}^-$, we propose to use the penalty functions method.

### 5.1 Unconstrained game with penalized cost functional

Let us start with the $\text{ERG}^+$. For this game, the set of admissible evader's controls is $U_e^+ = \left\{ u_e(\cdot) \in L_2[0, t_f] : \ w_0 + \int_0^{t_f} g_e(t) u_e(t) dt = \mu_e a_e^{\max} \right\}$.

Consider an auxiliary linear-quadratic differential game for the system (24) – (25) with the cost functional

$$J_\varepsilon^+ = J - \frac{1}{\varepsilon}(w(t_f) - \mu_e a_e^{\max})^2, \tag{55}$$

where $\varepsilon > 0$ is a parameter. This game is called the Unconstrained Penalized Game corresponding to the $\text{ERG}^+$ ($\text{UPG}_\varepsilon^+$). Note that the $\text{UPG}_\varepsilon^+$, having no state constraints, is a standard linear-quadratic differential game. The Hamiltonian in this game is

$$H = H(t, \lambda_z, \lambda_w, u_p, u_e) =$$
$$\lambda_z(h_p(t) u_p + h_e(t) u_e) + \lambda_w g_e(t) u_e + \alpha u_p^2 - \beta u_e^2, \tag{56}$$

where $\lambda_z$ and $\lambda_w$ are the co-state variables, satisfying the differential equations

$$\dot{\lambda}_z = -\frac{\partial H}{\partial z} = 0, \ \lambda_z(t_f) = 2z(t_f), \tag{57}$$

$$\dot{\lambda}_w = -\frac{\partial H}{\partial w} = 0, \quad \lambda_w(t_f) = -\frac{2}{\varepsilon}(w(t_f) - \mu_e a_e^{\max}). \tag{58}$$

The solutions of (57) – (58) are

$$\lambda_z^\varepsilon \equiv 2z_{f\varepsilon}^+, \quad \lambda_w^\varepsilon \equiv -2v_{f\varepsilon}^+, \tag{59}$$

where

$$z_{f\varepsilon}^+ \triangleq z(t_f), \quad v_{f\varepsilon}^+ \triangleq \frac{w(t_f) - \mu_e a_e^{\max}}{\varepsilon}. \tag{60}$$

Due to (Bryson & Ho, 1975) and by using (59), the UPG$_\varepsilon^+$ candidate open-loop saddle point is

$$u_{p\varepsilon}^+(t) = \arg\min_{u_p} H(t, \lambda_z, \lambda_w, u_p, u_e) = -\frac{1}{\alpha} h_p(t) z_{f\varepsilon}^+, \tag{61}$$

$$u_{e\varepsilon}^+(t) = \arg\max_{u_e} H(t, \lambda_z, \lambda_w, u_p, u_e) = \frac{1}{\beta} \left[ h_e(t) z_{f\varepsilon}^+ - g_e(t) v_{f\varepsilon}^+ \right]. \tag{62}$$

Substituting (61) – (62) into the system (24) – (25), solving the resulting system and taking into account the notation (60) yield

$$z_{f\varepsilon}^+ = z_0 - \left( \frac{1}{\alpha} \int_0^{t_f} h_p^2(t) dt - \frac{1}{\beta} \int_0^{t_f} h_e^2(t) dt \right) z_{f\varepsilon}^+ - \frac{1}{\beta} \left( \int_0^{t_f} h_e(t) g_e(t) dt \right) v_{f\varepsilon}^+,$$

$$\varepsilon v_{f\varepsilon}^+ + \mu_e a_e^{\max} = w_0 + \frac{1}{\beta} \left( \int_0^{t_f} h_e(t) g_e(t) dt \right) z_{f\varepsilon}^+ - \frac{1}{\beta} \left( \int_0^{t_f} g_e^2(t) dt \right) v_{f\varepsilon}^+.$$

These equations lead to the linear system with respect to $\omega_\varepsilon^+ \triangleq (z_{f\varepsilon}^+, v_{f\varepsilon}^+)^T$:

$$(G + D_\varepsilon) \omega_\varepsilon^+ = b^+, \tag{63}$$

where

$$G \triangleq \begin{bmatrix} G_1 & G_2 \\ -G_2 & G_3 \end{bmatrix}, \tag{64}$$

$$G_1 = s, \quad G_2 = \frac{1}{\beta} \int_0^{t_f} h_e(t) g_e(t) dt, \quad G_3 = \frac{1}{\beta} \int_0^{t_f} g_e^2(t) dt, \tag{65}$$

$$D_\varepsilon \triangleq \begin{bmatrix} 0 & 0 \\ 0 & \varepsilon \end{bmatrix}, \quad b^+ = \begin{bmatrix} z_0 \\ w_0 - \mu_e a_e^{\max} \end{bmatrix}, \tag{66}$$

$s$ is given by (32). It can be seen directly that, subject to the condition (33), the matrix $G + D_\varepsilon$ is non-singular for all $\varepsilon \geq 0$. Thus, the system (63) has the unique solution

$$\omega_\varepsilon^+ = (G + D_\varepsilon)^{-1} b^+. \tag{67}$$

Thus, the candidate optimal players' controls in the $\text{UPG}_\varepsilon^+$ are completely derived.

**Proposition 2** *If the condition (33) holds, then for any $\varepsilon > 0$, the pair (61) – (62) constitutes the $\text{UPG}_\varepsilon^+$ saddle point, i.e.,*

$$J_\varepsilon^+(u_{p\varepsilon}^+(\cdot), u_e(\cdot)) \leq J_\varepsilon^+(u_{p\varepsilon}^+(\cdot), u_{e\varepsilon}^+(\cdot)) \leq J_\varepsilon^+(u_p(\cdot), u_{e\varepsilon}^+(\cdot)), \tag{68}$$

*for all $u_p(\cdot), u_e(\cdot) \in L_2[0, t_f]$ and for $\varepsilon > 0$.*

**Proof.** For the $\text{UPG}_\varepsilon^+$, the Riccati equation (Bernhard, 2015, Eq. (2)) becomes

$$\dot{P} = -\frac{1}{\beta} P \begin{bmatrix} h_e^2(t) & h_e(t)g_e(t) \\ h_e(t)g_e(t) & g_e^2(t) \end{bmatrix} P, \quad P(t_f) = \begin{bmatrix} 1 & 0 \\ 0 & -1/\varepsilon \end{bmatrix}. \tag{69}$$

It is verified directly that the solution of (69) is

$$P(t) = \begin{bmatrix} 1 - \frac{1}{\beta} \int_t^{t_f} h_e^2(\xi) d\xi & -\frac{1}{\beta} \int_t^{t_f} h_e(\xi) g_e(\xi) d\xi \\ \frac{1}{\beta\varepsilon} \int_t^{t_f} h_e(\xi) g_e(\xi) d\xi & 1 + \frac{1}{\beta\varepsilon} \int_t^{t_f} g_e^2(\xi) d\xi \end{bmatrix}^{-1} \begin{bmatrix} 1 & 0 \\ 0 & -1/\varepsilon \end{bmatrix}. \tag{70}$$

The condition (33) guarantees the existence of the inverse matrix in (70) for all $t \in [0, t_f]$ and all $\varepsilon \geq 0$. Thus, by virtue of (Bernhard, 2015), the pair (61) – (62) indeed constitutes the $\text{UPG}_\varepsilon^+$ saddle point. $\square$

By substituting the saddle point (61) – (62) into the cost functional (55), the value of the $\text{UPG}_\varepsilon^+$ is

$$(J_\varepsilon^+)^* \triangleq J_\varepsilon^+(u_{p\varepsilon}^+(\cdot), u_{e\varepsilon}^+(\cdot)) = (\omega_\varepsilon^+)^T \tilde{G}_\varepsilon \omega_\varepsilon^+, \tag{71}$$

where $\tilde{G}_\varepsilon \triangleq \text{diag}(1, -1)(G + D_\varepsilon)$.

**Remark 5** For the $ERG^-$, the set of admissible evader's controls is

$$U_e^- = \left\{ u_e(\cdot) \in L_2[0, t_f] : \ w_0 + \int_0^{t_f} g_e(t) u_e(t) dt = -\mu_e a_e^{\max} \right\}.$$

For the $UPG_\varepsilon^-$, corresponding to the $ERG^-$, the cost functional is

$$J_\varepsilon^- = J - \frac{1}{\varepsilon}(w(t_f) + \mu_e a_e^{\max})^2.$$

The $UPG_\varepsilon^-$ saddle point $(u_{p\varepsilon}^-(\cdot), u_{e\varepsilon}^-(\cdot))$ has the same form as (61) – (62) where the vector $\omega_\varepsilon^+ = (z_{f\varepsilon}^+, v_{f\varepsilon}^+)^T$ is replaced by $\omega_\varepsilon^- = (z_{f\varepsilon}^-, v_{f\varepsilon}^-)^T$ satisfying the system

$$(G + D_\varepsilon) \omega_\varepsilon^- = b^-, \quad b^- = \begin{bmatrix} z_0 \\ w_0 + \mu_e a_e^{\max} \end{bmatrix}.$$

The value of the $UPG_\varepsilon^-$ is $(J_\varepsilon^-)^* = (\omega_\varepsilon^-)^T \tilde{G}_\varepsilon \omega_\varepsilon^-$.

## 5.2 Limit of the $UPG_\varepsilon^+/UPG_\varepsilon^-$ solutions for $\varepsilon \to 0$

**Lemma 3** If the condition (33) holds, then for the solution of the $UPG_\varepsilon^+$,

$$u_p^+(t) \triangleq \lim_{\varepsilon \to 0} u_{p\varepsilon}^+(t) = -\frac{1}{\alpha} h_p(t) z_f^+, \quad t \in [0, t_f], \tag{72}$$

$$u_e^+(t) \triangleq \lim_{\varepsilon \to 0} u_{e\varepsilon}^+(t) = \frac{1}{\beta} \left[ h_e(t) z_f^+ - g_e(t) v_f^+ \right], \quad t \in [0, t_f], \tag{73}$$

$$(J^+)^* \triangleq \lim_{\varepsilon \to 0} (J_\varepsilon^+)^* = (\omega_f^+)^T \tilde{G} \omega_f^+, \tag{74}$$

where

$$\omega_f^+ = (z_f^+, v_f^+)^T = G^{-1} b^+, \tag{75}$$

$$\tilde{G} = \text{diag}(1, -1) G. \tag{76}$$

**Proof.** By limiting (67) for $\varepsilon \to 0$,

$$\omega_f^+ \triangleq \lim_{\varepsilon \to 0} \omega_\varepsilon^+ = G^{-1} b^+. \tag{77}$$

Now, calculating the limits of the $UPG_\varepsilon^+$ saddle point (61) – (62) and the game value (71) for $\varepsilon \to 0$, and using (77) directly yield (72) – (74). □

**Corollary 1** *The value $(J^+)^*$ can be expressed as*

$$(J^+)^* = \chi_0^T \bar{G} \chi_0 + 2\chi_0^T \bar{G} \gamma^+ + (\gamma^+)^T \bar{G} \gamma^+, \tag{78}$$

*where*

$$\bar{G} \triangleq (G^{-1})^T \mathrm{diag}(1, -1), \tag{79}$$

$$\chi_0 \triangleq (z_0, w_0)^T, \tag{80}$$

$$\gamma^+ \triangleq (0, -\mu_e a_e^{\max})^T, \tag{81}$$

**Proof.** The expression (78) is proved by direct substitution of $\omega_f^+ = G^{-1} b^+$ into (74) and by the representation

$$b^+ = \chi_0 + \gamma^+. \tag{82}$$

$\square$

**Remark 6** *Similarly to Lemma 3, for the solution of the $UPG_\varepsilon^-$,*

$$u_p^-(t) \triangleq \lim_{\varepsilon \to 0} u_{p\varepsilon}^-(t) = -\frac{1}{\alpha} h_p(t) z_f^-, \quad t \in [0, t_f], \tag{83}$$

$$u_e^-(t) \triangleq \lim_{\varepsilon \to 0} u_{e\varepsilon}^-(t) = \frac{1}{\beta} \left[ h_e(t) z_f^- - g_e(t) v_f^- \right], \quad t \in [0, t_f], \tag{84}$$

$$(J^-)^* \triangleq \lim_{\varepsilon \to 0} (J_\varepsilon^-)^* = (\omega_f^-)^T \tilde{G} \omega_f^- = \chi_0^T \bar{G} \chi_0 + 2\chi_0^T \bar{G} \gamma^- + (\gamma^-)^T \bar{G} \gamma^-, \tag{85}$$

*where*

$$\omega_f^- = (z_f^-, v_f^-)^T = G^{-1} b^-, \quad \gamma^- \triangleq (0, \mu_e a_e^{\max})^T. \tag{86}$$

## 5.3 ERG$^+$/ERG$^-$ saddle points

**Theorem 1** *Let the condition (33) hold. Then,*

$$u_e^+(\cdot) \in U_e^+, \tag{87}$$

*i.e., for $u_e = u_e^+$, the equality constraint $w(t_f) = \mu_e a_e^{\max}$ is satisfied in the $ERG^+$. Moreover, the pair $(u_p^+(\cdot), u_e^+(\cdot))$ given by (72) – (73) constitutes the saddle point of the $ERG^+$, whereas the value $(J^+)^*$ given by (74) is the value of the $ERG^+$:*

$$J(u_p^+(\cdot), u_e(\cdot)) \leq J(u_p^+(\cdot), u_e^+(\cdot)) \leq J(u_p(\cdot), u_e^+(\cdot)), \quad \forall\, u_p(\cdot) \in U_p,\ u_e(\cdot) \in U_e^+, \tag{88}$$

*and*

$$J(u_p^+(\cdot), u_e^+(\cdot)) = (J^+)^*. \tag{89}$$

**Proof.** By integrating the differential equation (25) for $u_e = u_e^+$,

$$w(t_f) = w_0 + G_2 z_f^+ - G_3 v_f^+. \tag{90}$$

From the form of the matrix $G$, defined by (64) – (65), vector $b^+$, given in (66), and from (75),

$$\mu_e a_e^{\max} = w_0 + G_2 z_f^+ - G_3 v_f^+. \tag{91}$$

Comparison of (90) and (91) leads to $w(t_f) = \mu_e a_e^{\max}$, which implies the inclusion (87). Similarly, by integrating the differential equation (24) for $u_p = u_p^+$, $u_e = u_e^+$,

$$z(t_f) = z_f^+. \tag{92}$$

Substituting $u_p^+(\cdot)$ and $u_e^+(\cdot)$ into the original cost functional $J$ (see (27)) and using (92) yield after a routine algebra the equality (89).

Now, let us proceed to the proof of the ERG$^+$ saddle point inequality (88). We prove it by limiting for $\varepsilon \to 0$ in the UPG$_\varepsilon^+$ saddle inequality (68). Due to (24), (55), (61), for $u_p(\cdot) = u_{p\varepsilon}^+(\cdot)$ and any $u_e(\cdot) \in U_e^+ \subset U_e$,

$$J_\varepsilon^+(u_{p\varepsilon}^+(\cdot), u_e(\cdot)) = \tilde{z}^2(t_f) + \frac{1}{\alpha} \int_0^{t_f} h_p^2(t) dt (z_f^\varepsilon)^2 - \beta \int_0^{t_f} u_e^2(t) dt, \tag{93}$$

where $\tilde{z}(t)$ is the solution of the differential equation (24) generated by $u_{p\varepsilon}^+(\cdot)$ and $u_e(\cdot)$. Thus,

$$\tilde{z}(t_f) = z_0 - \frac{1}{\alpha} \int_0^{t_f} h_p^2(t) dt (z_f^\varepsilon) + \int_0^{t_f} h_e(t) u_e(t) dt. \tag{94}$$

Due to (75) and (77),

$$\lim_{\varepsilon \to 0} z_f^\varepsilon = z_f^+, \quad \lim_{\varepsilon \to 0} v_f^\varepsilon = v_f^+. \tag{95}$$

Using (27), (72) and (93) – (95) yields

$$\lim_{\varepsilon \to 0} J_\varepsilon^+(u_{p\varepsilon}^+(\cdot), u_e(\cdot)) = J(u_p^+(\cdot), u_e(\cdot)), \quad u_e(\cdot) \in U_e^+ \tag{96}$$

Similarly it is shown that

$$\lim_{\varepsilon \to 0} J_\varepsilon^+(u_p(\cdot), u_{e\varepsilon}^+(\cdot)) = J(u_p(\cdot), u_e^+(\cdot)), \quad u_p(\cdot) \in U_p. \tag{97}$$

Furthermore, by (71), (74) and (89),

$$\lim_{\varepsilon \to 0} J_\varepsilon^+(u_{p\varepsilon}^+(\cdot), u_{e\varepsilon}^+(\cdot)) = J(u_p^+(\cdot), u_e^+(\cdot)). \tag{98}$$

Now, based on (96) – (98), the limiting $\varepsilon \to 0$ in (68) for $u_p(\cdot) \in U_p$, $u_e(\cdot) \in U_e^+$ leads to the ERG$^+$ saddle-point inequality (88). This completes the proof of the theorem. $\square$

**Remark 7** *Similarly to Theorem 1, subject to the condition (33), $u_e^-(\cdot) \in U_e^-$, i.e., for $u_e = u_e^-$, the equality constraint $w(t_f) = -\mu_e a_e^{\max}$ is satisfied in the ERG$^-$. Moreover, the pair $(u_p^-(\cdot), u_e^-(\cdot))$ given by (83) – (84) constitutes the saddle point of the ERG$^-$, whereas the value $(J^-)^*$ given by (85) is the value of the ERG$^-$:*

$$J(u_p^-(\cdot), u_e(\cdot)) \leq J(u_p^-(\cdot), u_e^-(\cdot)) \leq J(u_p(\cdot), u_e^-(\cdot)), \quad \forall \, u_p(\cdot) \in U_p, \, u_e(\cdot) \in U_e^-,$$

*and $J(u_p^-(\cdot), u_e^-(\cdot)) = (J^-)^*$. Moreover, in the ERG$^-$,*

$$w(t_f) = -\mu a_e^{\max}, \tag{99}$$

$$z(t_f) = z_f^-. \tag{100}$$

# 6 Saddle Point Solution of the Equality Constraint Game for $(z_0, w_0) \notin \Omega$

The set (49) of admissible evader's controls in the ERG can be represented as

$$\bar{U}_e = U_e^+ \cup U_e^-. \tag{101}$$

In this section, we establish the conditions guaranteeing that the ERG saddle point coincides with the saddle point either of the ERG$^+$, or of the ERG$^-$. This means that one of the saddle point inequalities is valid:

$$J(u_p^+(\cdot), u_e(\cdot)) \leq J(u_p^+(\cdot), u_e^+(\cdot)) \leq J(u_p(\cdot), u_e^+(\cdot)), \quad \forall \, u_p(\cdot) \in U_p, \, u_e(\cdot) \in \bar{U}_e, \tag{102}$$

or

$$J(u_p^-(\cdot), u_e(\cdot)) \leq J(u_p^-(\cdot), u_e^-(\cdot)) \leq J(u_p(\cdot), u_e^-(\cdot)), \quad \forall \, u_p(\cdot) \in U_p, \, u_e(\cdot) \in \bar{U}_e, \tag{103}$$

The derivation of these conditions is based on two auxiliary optimal control problems with terminal state equality constraints. The first problem is formulated for the system

$$\dot{z} = h_p(t)u_p^+(t) + h_e(t)u_e, \quad z(0) = z_0, \tag{104}$$

$$\dot{w} = g_e(t)u_e, \quad w(0) = w_0, \quad w(t_f) = -\mu_e a_e^{\max}, \tag{105}$$

and the cost functional

$$J_{e1}(u_p^+(\cdot), u_e(\cdot)) \triangleq z^2(t_f) - \beta \int_0^{t_f} u_e^2(t)dt \to \max_{u_e(\cdot) \in U_e^-}. \tag{106}$$

The second problem is formulated for the system and the cost functional

$$\dot{z} = h_p(t)u_p^-(t) + h_e(t)u_e, \quad z(0) = z_0, \tag{107}$$

$$\dot{w} = g_e(t)u_e, \quad w(0) = w_0, \quad w(t_f) = \mu_e a_e^{\max}, \tag{108}$$

$$J_{e2}(u_p^-(\cdot), u_e(\cdot)) \triangleq z^2(t_f) - \beta \int_0^{t_f} u_e^2(t)dt \to \max_{u_e(\cdot) \in U_e^+}. \tag{109}$$

## 6.1 Solution of the auxiliary optimal control problems with terminal state equality constraints

Due to the results of (Ioffe & Tikhomirov, 1979), both auxiliary optimal control problems have solutions.

### 6.1.1 Solution of the first problem

Let us start with the problem (104) – (106). By the Pontryagin Maximum Principle (Vasilyev, 1988), the Hamiltonian of this problem is

$$H_1 = \lambda_z(h_p(t)u_p^+(t) + h_e(t)u_e) + \lambda_w g_e(t)u_e - \beta u_e^2,$$

where the co-states $\lambda_z$ and $\lambda_w$ satisfy

$$\dot{\lambda}_z = -\frac{\partial H}{\partial z} = 0, \quad \lambda_z(t_f) = 2z(t_f); \quad \dot{\lambda}_w = -\frac{\partial H}{\partial w} = 0.$$

Thus, $\lambda_z \equiv 2z(t_f)$, $\lambda_w \equiv \text{const}$.

The optimal control is

$$u_{e1}^*(t) = \arg\max_{u_e \in U_e^-} H = \frac{h_e(t)\lambda_z + g_e(t)\lambda_w}{2\beta} = \frac{1}{\beta}(h_e(t)z(t_f) - g_e(t)\tilde{\lambda}_w), \quad (110)$$

where $\tilde{\lambda}_w \triangleq -\lambda_w/2$.

Substituting (110) into the equations (104) – (105), following by their integration, yield the linear algebraic system w.r.t. $\omega_1 \triangleq (z(t_f), \tilde{\lambda}_w)^T$:

$$F\omega_1 = \mu_1, \quad (111)$$

where

$$F \triangleq \begin{bmatrix} 1 - \dfrac{1}{\beta}\displaystyle\int_0^{t_f} h_e^2(t)dt & G_2 \\ -G_2 & G_3 \end{bmatrix}, \quad \mu_1 \triangleq \chi_0 + \xi_1, \quad \xi_1 \triangleq \gamma^- + \begin{bmatrix} \displaystyle\int_0^{t_f} h_p(t)u_p^+(t)dt \\ 0 \end{bmatrix}. \quad (112)$$

Remember that the values $G_2$ and $G_3$ are defined in (65), the vectors $\chi_0$ and $\gamma^-$ are given by (80) and in (86), respectively.

Subject to the condition (33), the matrix $F$ is non-singular, yielding the unique solution of the system (111)

$$\omega_1 = F^{-1}\mu_1. \quad (113)$$

Substituting (110) into the functional (106) and taking into account (113) yield, after a routine algebra, the optimal value of the cost functional in the first auxiliary optimal control problem:

$$J_{e1}^* = J_{e1}(u_p^+(\cdot), u_{e1}^*(\cdot)) = \chi_0^T \bar{F} \chi_0 + 2\chi_0^T \bar{F} \xi_1 + \xi_1^T \bar{F} \xi_1, \quad (114)$$

where

$$\bar{F} = (F^{-1})^T \text{diag}(1, -1). \quad (115)$$

Due to (27) and (114),

$$J(u_p^+(\cdot), u_{e1}^*(\cdot)) = J_{e1}(u_p^+(\cdot), u_{e1}^*(\cdot)) + \alpha \int_0^{t_f} (u_p^+(t))^2 dt$$

$$= \chi_0^T \bar{F} \chi_0 + 2\chi_0^T \bar{F} \xi_1 + \xi_1^T \bar{F} \xi_1 + \alpha \int_0^{t_f} (u_p^+(t))^2 dt. \quad (116)$$

Let us show that the value (116) is a quadratic form w.r.t. the vector $\chi_0$. Indeed, due to (72), (75) and (82),

$$u_p^+(t) = -\frac{1}{\alpha} h_p(t)(1,0) G^{-1}(\chi_0 + \gamma^+). \tag{117}$$

By using (117), we can represent the vector $\xi_1$ given in (112) as

$$\xi_1 = \gamma^- - \mathrm{diag}(\nu_p, 0) G^{-1}(\chi_0 + \gamma^+), \tag{118}$$

where

$$\nu_p \triangleq \frac{1}{\alpha} \int_0^{t_f} h_p^2(t) dt. \tag{119}$$

Moreover, the last term of the right-hand side in (116) can be expressed as

$$\alpha \int_0^{t_f} (u_p^+(t))^2 dt = (\chi_0 + \gamma^+)^T (G^{-1})^T \mathrm{diag}(\nu_p, 0) G^{-1}(\chi_0 + \gamma^+). \tag{120}$$

Substitution of the equations (118) and (120) into (116) converts the latter to a quadratic form with respect to $\chi_0$:

$$J(u_p^+(\cdot), u_{e1}^*(\cdot)) = \chi_0^T \bar{G} \chi_0 + 2\chi_0^T \bar{G} \gamma^- + \rho, \tag{121}$$

where the matrix $\bar{G}$, given by (79), defines the quadratic term in the expression (78) for the ERG$^+$ value $(J^+)^*$,

$$\rho = \frac{(\mu_e a_e^{\max})^2 \left[3\nu_p G_2^2 - (G_1 - \nu_p) \det(G)\right]}{\det(G) \det(F)}.$$

### 6.1.2 Solution of the second problem

Similarly to the solution of the problem (104) – (106) we obtain the solution of the problem (107) – (109). Namely, the optimal control in this problem has the form $u_{e2}^*(t) = \frac{1}{\beta}(h_e(t), -g_e(t)) F^{-1} \mu_2$, where the matrix $F$ is given in (112), and $\mu_2 \triangleq \chi_0 + \xi_2$, $\xi_2 \triangleq \gamma^+ + \begin{bmatrix} \int_0^{t_f} h_p(t) u_p^-(t) dt \\ 0 \end{bmatrix}$.

The optimal value of the cost functional in the problem (107) – (109) is $J_{e2}^* = J_{e2}(u_p^-(\cdot), u_{e2}^*(\cdot)) = \chi_0^T \bar{F} \chi_0 + 2\chi_0^T \bar{F} \xi_2 + \xi_2^T \bar{F} \xi_2$, where the matrix $\bar{F}$ is given by (115).

The value of the cost functional (27), calculated for $u_p(t) = u_p^-(t)$ and $u_e(t) = u_{e2}^*(t)$, has the form

$$J(u_p^-(\cdot), u_{e2}^*(\cdot)) = J_{e2}(u_p^-(\cdot), u_{e2}^*(\cdot)) + \alpha \int_0^{t_f} (u_p^-(t))^2 dt$$

$$= \chi_0^T \bar{F} \chi_0 + 2\chi_0^T \bar{F} \xi_2 + \xi_2^T \bar{F} \xi_2 + \alpha \int_0^{t_f} (u_p^-(t))^2 dt. \tag{122}$$

Moreover, similarly to (117) – (120), we obtain $u_p^-(t) = -\frac{1}{\alpha} h_p(t)(1,0) G^{-1}(\chi_0 + \gamma^-)$,

$$\xi_2 = \gamma^+ - \text{diag}(\nu_p, 0) G^{-1}(\chi_0 + \gamma^-), \tag{123}$$

$$\alpha \int_0^{t_f} (u_p^-(t))^2 dt = (\chi_0 + \gamma^-)^T (G^{-1})^T \text{diag}(\nu_p, 0) G^{-1}(\chi_0 + \gamma^-). \tag{124}$$

Substitution of the equations (123) and (124) into (122) converts the latter to a quadratic form with respect to $\chi_0$:

$$J(u_p^-(\cdot), u_{e2}^*(\cdot)) = \chi_0^T \bar{G} \chi_0 + 2\chi_0^T \bar{G} \gamma^+ + \rho. \tag{125}$$

## 6.2 Conditions for fulfilment of the saddle point inequalities (102) and (103)

**Theorem 2** *Let the condition (33) hold. Let $(z_0, w_0) \notin \Omega$. Then, the pair $(u_p^+(\cdot), u_e^+(\cdot))$ constitutes a saddle point in the ERG if and only if*

$$J(u_p^+(\cdot), u_{e1}^*(\cdot)) \leq J(u_p^+(\cdot), u_e^+(\cdot)). \tag{126}$$

**Proof.** *Necessity.* Let the pair $(u_p^+(\cdot), u_e^+(\cdot))$ be a saddle point in the ERG. Then, the inequality (102) is valid. Due to the left-hand inequality in (102) and the equations (27), (106), $u_e^+(t)$ is an optimal control in the problem consisting of the equations of dynamics (104) and

$$\dot{w} = g_e(t) u_e, \quad w(0) = w_0, \quad |w(t_f)| = \mu_e a_e^{\max}, \tag{127}$$

and the performance index

$$J_{e1}(u_p^+(\cdot), u_e(\cdot)) \to \max_{u_e(\cdot) \in \bar{U}_e}. \qquad (128)$$

Remember, that $u_{e1}^*(t)$ is an optimal control in the problem (104) – (106). Comparing this optimal control problem and the problem (104), (127), (128), and taking into account the inclusion $U_e^- \subset \bar{U}_e$ (see (101)), we obtain the inequality $J_{e1}(u_p^+(\cdot), u_{e1}^*(\cdot)) \leq J_{e1}(u_p^+(\cdot), u_e^+(\cdot))$. The latter, along with the equations (27), (106), directly yields the inequality (126). This completes the proof of the necessity.

*Sufficiency.* Since $(u_p^+(\cdot), u_e^+(\cdot))$ is a saddle point in the ERG$^+$, then the inequality (102) is fulfilled for all $u_e(\cdot) \in U_e^+$. Let us show that this inequality is fulfilled for all $u_e(\cdot) \in U_e^-$. Note, that the right-hand side inequality in (102) is independent on $u_e(\cdot)$, and it is fulfilled for all $u_p(\cdot) \in U_p$. Thus, to prove the fulfilment of (102) for all $u_e(\cdot) \in U_e^-$, it is sufficient to show the validity of its left-hand inequality for these $u_e(\cdot)$.

Since $u_{e1}^*(\cdot)$ is the optimal control in the first auxiliary problem (104) – (106), then

$$J_{e1}(u_p^+(\cdot), u_e(\cdot)) \leq J_{e1}(u_p^+(\cdot), u_{e1}^*(\cdot)) \quad \forall \ u_e(\cdot) \in U_e^-. \qquad (129)$$

From (27) and (106), we directly obtain

$$J(u_p^+(\cdot), u_e(\cdot)) = J_{e1}(u_p^+(\cdot), u_e(\cdot)) + \alpha \int_0^{t_f} (u_p^+(t))^2 dt \quad \forall u_e(\cdot) \in U_e^-. \qquad (130)$$

The equations (116), (130) and the inequalities (126), (129) immediately yield the validity of the left-hand inequality in (102) for all $u_e(\cdot) \in U_e^-$. This completes the proof of the sufficiency. Thus, the theorem is proven. $\square$

**Theorem 3** *Let the condition (33) hold. Let $(z_0, w_0) \notin \Omega$. Then, the pair $(u_p^-(\cdot), u_e^-(\cdot))$ constitutes a saddle point in the ERG if and only if*

$$J(u_p^-(\cdot), u_{e2}^*(\cdot)) \leq J(u_p^-(\cdot), u_e^-(\cdot)). \qquad (131)$$

**Proof.** Using the solution of the second auxiliary problem (107) – (109) (see Section 6.1.2), the theorem is proven similarly to Theorem 2. $\square$

**Remark 8** *Due to (78) and (121), the condition (126) reads*

$$2\chi_0^T \bar{G}(\gamma^- - \gamma^+) + \rho - (\gamma^+)^T \bar{G}\gamma^+ \leq 0. \tag{132}$$

*By routine algebra, this inequality becomes*

$$w_0 + az_0 \geq d\mu_e a_e^{\max}, \tag{133}$$

*where the coefficient $a$ is the same as in the definition (37) of the set $\Omega$,*

$$d = \frac{\nu_p G_2^2}{G_1 \det(F)} > 0. \tag{134}$$

*The condition (133), along with the condition $(z_0, w_0) \notin \Omega$, define the set of initial positions*

$$\Omega^+ = \{(z_0, w_0) : |w_0 + az_0| > \mu_e a_e^{\max}, \ w_0 + az_0 \geq d\mu_e a_e^{\max}\}, \tag{135}$$

*for which, subject to (33), the pair $(u_p^+, u_e^+)$ is the saddle point in the ERG. Similarly, by using (85), (125), the condition (131), along with the condition $(z_0, w_0) \notin \Omega$, define the set of initial positions,*

$$\Omega^- = \{(z_0, w_0) : |w_0 + az_0| > \mu_e a_e^{\max}, \ w_0 + az_0 \leq -d\mu_e a_e^{\max}\}, \tag{136}$$

*for which, subject to (33), the pair $(u_p^-, u_e^-)$ is the saddle point in the ERG.*

**Proposition 3** *Subject to the condition (33),*

$$\Omega^+ = \{(z_0, w_0) : w_0 + az_0 \geq \mu_e a_e^{\max}\}, \tag{137}$$

$$\Omega^- = \{(z_0, w_0) : w_0 + az_0 \leq -\mu_e a_e^{\max}\}, \tag{138}$$

**Proof.** Due to (135) – (136), it is sufficient to prove that

$$d < 1. \tag{139}$$

By virtue of (134), this means that

$$G_1 \det F - \nu_p G_2^2 > 0. \tag{140}$$

Due to (112), $\det F = (G_1 - \nu_p)G_3 + G_2^2$, and (140) reads

$$(G_1 - \nu_p)(G_1 G_3 + G_2^2) = (G_1 - \nu_p) \det G > 0. \tag{141}$$

By (32), (33), (65) and (119),

$$G_1 - \nu_p = 1 - \frac{1}{\beta} \int_0^{t_f} h_e^2(t) dt > 0. \tag{142}$$

Moreover, it was shown that $\det G > 0$, which, along with (142), leads to (139). This completes the proof of the proposition

$\square$

# 7 Saddle Point Solution in the Inequality Constraint Game for $(z_0, w_0) \notin \Omega$

In this section, we show that the ERG saddle points $(u_p^+(\cdot), u_e^+(\cdot))$ and $(u_p^-(\cdot), u_e^-(\cdot))$, obtained in the previous section, also are the saddle points in the RG. To do this, we are based on Lemma 2. In accordance with this lemma, we consider two auxiliary optimal control problems with terminal state inequality constraints and for $(z_0, w_0) \notin \Omega$. These problems are obtained from the optimal control problem (50) – (52), (28) by replacing there the pursuer's control $\bar{u}_p^*(t)$ with the controls $u_p^+(t)$ and $u_p^-(t)$, respectively.

## 7.1 Solution of the auxiliary optimal control problems with terminal state inequality constraints

Due to the results of (Ioffe & Tikhomirov, 1979), these optimal control problems have solutions, the optimal controls $\bar{u}_{e1}^*(t)$ and $\bar{u}_{e2}^*(t)$, respectively.

### 7.1.1 Solution of the first problem

We start with the problem (50) – (52), (28) where $\bar{u}_p^*(t) = u_p^+(t)$, $t \in [0, t_f]$, and $(z_0, w_0) \notin \Omega$. Remember that in Theorem 2 the ERG saddle points $(u_p^+(\cdot), u_e^+(\cdot))$ was obtained subject to the inequality (126). Therefore, in this section we assume this inequality to be valid.

Similarly to the proof of Lemma 1 (see the equations (39) – (42)), there exist the functions $\bar{\lambda}_z(t)$, $\bar{\lambda}_w(t)$ and the numbers $\bar{\lambda}_1^f \geq 0$, $\bar{\lambda}_2^f \geq 0$, such that

$$\dot{\bar{\lambda}}_z = 0, \quad \bar{\lambda}_z(t_f) = 2z(t_f), \tag{143}$$

$$\dot{\bar{\lambda}}_w = 0, \quad \bar{\lambda}_w(t_f) = \bar{\lambda}_1^f - \bar{\lambda}_2^f, \tag{144}$$

$$\bar{\lambda}_1^f(w(t_f) - \mu_e a_e^{\max}) = 0, \quad \bar{\lambda}_2^f(w(t_f) + \mu_e a_e^{\max}) = 0, \tag{145}$$

and the optimal control satisfies the equation

$$\bar{u}_{e1}^*(t) = \frac{1}{2\beta}\left[h_e(t)\bar{\lambda}_z(t) + g_e(t)\bar{\lambda}_w(t)\right], \tag{146}$$

where $z(t)$ and $w(t)$ are the solutions of (50) – (51), generated by the pair $(u_p^+(\cdot), \bar{u}_{e1}^*(\cdot))$. Note that (143) yields

$$\bar{\lambda}_z(t) \equiv 2z(t_f). \tag{147}$$

For the further analysis of the set (144) – (147), the following three cases should be distinguished: (i) the optimal control $\bar{u}_{e1}^*(t)$ provides the fulfilment of the equality $w(t_f) = \mu_e a_e^{\max}$; (ii) the optimal control $\bar{u}_{e1}^*(t)$ provides the fulfilment of the equality $w(t_f) = -\mu_e a_e^{\max}$; (iii) the optimal control $\bar{u}_{e1}^*(t)$ provides the fulfilment of the inequality $-\mu_e a_e^{\max} < w(t_f) < \mu_e a_e^{\max}$.

**Case (i).** Since the pair $(u_p^+(\cdot), u_e^+(\cdot))$ is the saddle point in the ERG$^+$, then in this case the optimal control is $\bar{u}_{e1}^*(t) = u_e^+(t)$, $t \in [0, t_f]$. The value of the cost functional for this control is $\bar{J}_e(u_p^+(\cdot), u_e^+(\cdot)) = J(u_p^+(\cdot), u_e^+(\cdot)) - \alpha \int_0^{t_f} \left(u_p^+(t)\right)^2 dt.$

**Case (ii).** Since $u_{e1}^*(t)$ is the optimal control of the problem (104) – (106) (see Section 6.1.1), then in this case the optimal control is $\bar{u}_{e1}^*(t) = u_{e1}^*(t)$, $t \in [0, t_f]$. The value of the cost functional for this control is $\bar{J}_e(u_p^+(\cdot), u_{e1}^*(\cdot)) = J_{e1}(u_p^+(\cdot), u_{e1}^*(\cdot)) = J(u_p^+(\cdot), u_{e1}^*(\cdot)) - \alpha \int_0^{t_f} \left(u_p^+(t)\right)^2 dt.$

**Case (iii).** In this case, due to (145) and (144), we obtain $\bar{\lambda}_w(t) \equiv 0$. Therefore, by virtue of (146)–(147),

$$\bar{u}_{e1}^*(t) = \frac{1}{\beta}h_e(t)z(t_f). \tag{148}$$

Substitution of $u_p^+(t)$ (see (117)) and $\bar{u}_{e1}^*(t)$ into (50) instead of $\bar{u}_p^*(t)$ and $\bar{u}_e$, respectively, and the integration of the resulting initial-value problem yield the equation with respect to $z(t_f)$

$$z(t_f) = z_0 - (\nu_p, 0)G^{-1}(\chi_0 + \gamma^+) + \nu_e z(t_f), \quad \nu_e \triangleq \frac{1}{\beta} \int_0^{t_f} h_e^2(t)dt.$$

Subject to the condition (33), this equation has the unique solution

$$z(t_f) = \bar{z}_{f1} \triangleq \frac{1}{1 - \nu_e} \left( z_0 - (\nu_p, 0)G^{-1}(\chi_0 + \gamma^+) \right).$$

Thus, the control (148) becomes $\bar{u}_{e1}^*(t) = \frac{1}{\beta} h_e(t)\bar{z}_{f1}$. Substitution of this control into (51) instead of $\bar{u}_e$ and integration of the resulting initial-value problem yields $w(t_f) = w_0 + G_2 \bar{z}_{f1}$, where the value $G_2$ is given in (65). Using this expression for $w(t_f)$, one can directly conclude that in Case (iii) the following inequality is fulfilled:

$$-\mu_e a_e^{\max} < w_0 + G_2 \bar{z}_{f1} < \mu_e a_e^{\max}. \tag{149}$$

By routine algebra, the inequality (149) is rewritten as

$$(-1 + 2d)\mu_e a_e^{\max} < w_0 + a z_0 < \mu_e a_e^{\max}, \tag{150}$$

where $d$ is defined by (134). By (139), $d$ satisfies $0 < d < 1$. Then, $-1 < -1 + 2d < 1$, meaning that positions $(z_0, w_0)$ satisfying (149), belong to $\Omega$. Thus, for $(z_0, w_0) \notin \Omega$, Case (iii) is non-feasible.

Comparing the analysis results of Cases (i) – (iii) and using the equation (27) and the inequality (126), we directly obtain the assertion.

**Proposition 4** *Let the condition (33) and the inequality (126) hold. Let $(z_0, w_0) \notin \Omega$. Then Case (i) is valid and $u_e^+(t)$ is the optimal control in the first auxiliary problem.*

### 7.1.2 Solution of the second problem

Now, we consider the problem (50) – (52), (28) where $\bar{u}_p^*(t) = u_p^-(t)$, $t \in [0, t_f]$, and $(z_0, w_0) \notin \Omega$. This problem is solved similarly to the optimal control problem of the previous section. We solve it subject to the inequality (131). Due to Theorem 3, this inequality provides the pair $(u_p^-(\cdot), u_e^-(\cdot))$ to be a saddle point in the ERG.

To formulate the assertion, similar to Proposition 4, we use the following values of the cost functional $\bar{J}(u_p^-(\cdot), \bar{u}_e(\cdot))$ in the second auxiliary control problem:

$$\bar{J}(u_p^-(\cdot), \bar{u}_e^-(\cdot)) = J(u_p^-(\cdot), u_e^-(\cdot)) - \alpha \int_0^{t_f} \left(u_p^-(t)\right)^2 dt,$$

$$\bar{J}(u_p^-(\cdot), \bar{u}_{e2}^*(\cdot)) = \bar{z}_{f2}^2 - \frac{\bar{z}_{f2}}{\beta} \int_0^{t_f} h_e^2(t) dt,$$

where $\bar{u}_{e2}^*(t) = (1/\beta) h_e(t) \bar{z}_{f2}$, $\bar{z}_{f2} \triangleq \dfrac{1}{1-\nu_e} \left(z_0 - (\nu_p, 0) G^{-1}(\chi_0 + \gamma^-)\right)$.

Similarly to Proposition 4, we obtain the following proposition.

**Proposition 5** *Let the condition (33) and the inequality (131) hold. Let $(z_0, w_0) \notin \Omega$. Then $u_e^-(t)$ is the optimal control in the second auxiliary problem.*

## 7.2 Saddle points in the Reduced Game (24) − (25), (27), (28)

Based on Lemma 2, Theorems 2, 3 and Propositions 3 – 5, we immediately have the following two theorems.

**Theorem 4** *Let the condition (33) hold and $(z_0, w_0) \in \Omega^+$. Then, the pair $(u_p^+(\cdot), u_e^+(\cdot))$ is an open-loop saddle point in the Reduced Game.*

**Theorem 5** *Let the condition (33) hold and $(z_0, w_0) \in \Omega^-$. Then, the pair $(u_p^-(\cdot), u_e^-(\cdot))$ is an open-loop saddle point in the Reduced Game.*

**Remark 9** *Remember that in Section 3 the following result was obtained. If the condition (33) holds and $(z_0, w_0) \in \Omega$, then the pair $\left(u_p^0(\cdot), u_e^0(\cdot)\right)$ is an open-loop saddle point in the Reduced Game (24) – (25), (27), (28). Since the planar sets $\Omega$, $\Omega^+$ and $\Omega^-$ do not intersect each other, and $\Omega \cup \Omega^+ \cup \Omega^-$ coincides with the entire $(z_0, w_0)$-plane, then the above mentioned result of Section 3 and Theorems 4 – 5 represent the complete open-loop saddle-point solution of the Reduced Game.*

# 8 Special Case: First-Order Pursuer against First-Order Evader

In this section, the theory of the previous sections is applied to the particular case of (1) – (2) which is of a practical interest. This example illustrates some important features of the game solution.

## 8.1 Original Game OG

If both the pursuer and the evader have the first-order dynamics controller, then in the system (1) – (2), $n_p = 1$, $\bar{A}_p = -1/\tau_p$, $\bar{b}_p = 1/\tau_p$, $\bar{d}_p = 0$, $n_e = 1$, $\bar{A}_e = -1/\tau_e$, $\bar{b}_e = 1/\tau_e$, $\bar{d}_e = 0$, where $\tau_p$ and $\tau_e$ are the time constants of the pursuer's and the evader's controllers. The pursuer's and the evader's controls are the lateral acceleration commands.

In the OG, the controlled system is given by (7), where, by virtue of (4) – (6) and (8) – (9), $x = (y_p, \dot{y}_p, a_p, y_e, \dot{y}_e, a_e)^T$,

$$A = \begin{bmatrix} 0 & 1 & 0 & 0 & 0 & 0 \\ 0 & 0 & 1 & 0 & 0 & 0 \\ 0 & 0 & -1/\tau_p & 0 & 0 & 0 \\ 0 & 0 & 0 & 0 & 1 & 0 \\ 0 & 0 & 0 & 0 & 0 & 1 \\ 0 & 0 & 0 & 0 & 0 & -1/\tau_e \end{bmatrix}, \; B = \begin{bmatrix} 0 \\ 0 \\ 1/\tau_p \\ 0 \\ 0 \\ 0 \end{bmatrix}, \; C = \begin{bmatrix} 0 \\ 0 \\ 0 \\ 0 \\ 0 \\ 1/\tau_e \end{bmatrix}.$$

The cost functional (10) becomes

$$J = (x_4(t_f) - x_1(t_f))^2 + \alpha \int_0^{t_f} u_p^2(t)dt - \beta \int_0^{t_f} u_e^2(t)dt.$$

The matrix $\Phi_e$, given by (13) – (14), is

$$\Phi_e(t_f + t_c, t) = \begin{bmatrix} 1 & t_f + t_c - t & -\tau_e^2 \psi((t_f + t_c - t)/\tau_e) \\ 0 & 1 & \tau_e[\exp(-(t_f + t_c - t)/\tau_e) - 1] \\ 0 & 0 & \exp(-(t_f + t_c - t)/\tau_e) \end{bmatrix},$$

where $\psi(t) \triangleq \exp(-t) + t - 1 \geq 0$.

Thus, the terminal inequality constraint (16) becomes

$$|x_4(t_f) + t_c x_5(t_f) - \tau_e^2 \psi(t_c/\tau_e) x_6(t_f)| \leq \mu_e a_e^{\max},$$

where

$$\mu_e = \tau_e \int_{t_f}^{t_f+t_c} \psi((t_f + t_c - t)/\tau_e) dt =$$

$$\tau_e^2(1 - \sigma + \sigma^2/2 - \exp(-\sigma)), \quad \sigma \triangleq t_c/\tau_e. \tag{151}$$

## 8.2 Reduced Game

In this example, $D_{ep} = [1, 0, 0, 0]$ and the matrix $\Phi_{ep}$, introduced in (20), is

$$\Phi_{ep}(t, \tau) = \begin{bmatrix} 1 & t - \tau & -\tau_p^2 \psi((t-\tau)/\tau_p) & \tau_e^2 \psi((t-\tau)/\tau_e) \\ 0 & 1 & \tau_p[\exp(-(t-\tau)/\tau_p) - 1] & -\tau_e[\exp(-(t-\tau)/\tau_e) - 1] \\ 0 & 0 & \exp(-(t-\tau)/\tau_p) & 0 \\ 0 & 0 & 0 & \exp(-(t-\tau)/\tau_e) \end{bmatrix}.$$

Thus, the scalar variables (20) and (23) become

$$z(t) = y_e - y_p + (t_f - t)(\dot{y}_e - \dot{y}_p) - \tau_p^2 \psi((t_f - t)/\tau_p)a_p + \tau_e^2 \psi((t_f - t)/\tau_e)a_e,$$

$$w(t) = y_e + (t_f + t_c - t)\dot{y}_e + \tau_e^2 \psi((t_f + t_c - t)/\tau_e)a_e.$$

The coefficient functions (26) in the differential equations (24) – (25) become

$$h_p(t) = -\tau_p \psi((t_f - t)/\tau_p), \quad h_e(t) = \tau_e \psi((t_f - t)/\tau_e), \tag{152}$$

$$g_e(t) = \tau_e \psi((t_f + t_c - t)/\tau_e). \tag{153}$$

The differential equations (24) – (25) become

$$\dot{z} = -\tau_p \psi((t_f - t)/\tau_p)u_p + \tau_e \psi((t_f - t)/\tau_e)u_e, \tag{154}$$

$$\dot{w} = \tau_e \psi((t_f + t_c - t)/\tau_e)u_e. \tag{155}$$

The Reduced Game (RG) is formulated for the system (154) – (155) with the cost functional (27) and the terminal inequality constraint (28) where $\mu_e$ is given by (151).

## 8.3 Unconstrained Reduced Game

In the Unconstrained Reduced Game (URG), the optimal controls (34) are $u_p^0(t) = \dfrac{\tau_p \psi((t_f - t)/\tau_p) z_0}{\alpha s}$, $u_e^0(t) = \dfrac{\tau_e \psi((t_f - t)/\tau_e) z_0}{\beta s}$, where $s$ is defined by (32) by substituting $h_p(t)$ and $h_e(t)$ from (152). The solvability condition (33) reads

$$\beta > \beta^* = \tau_e^2 \int_0^{t_f} \psi^2((t_f - t)/\tau_e) dt, \tag{156}$$

yielding $s > 0$.

**Remark 10** *By comparing (36) and (65), the coefficient $a$ in the definition (37) of the set $\Omega$ writes as $a = -\dfrac{G_2}{G_1}$.*

To illustrate numerically the above presented results, we choose $t_f = 1$ s, $\nu = 0.9$, $a_e^{\max} = 100$ m/s$^2$, $\beta = 0.3$, $\tau_p = 0.2$, $\tau_e = 0.1$ s. For these parameters, $\beta^* = 0.2438$, and the solvability condition (156) is valid.

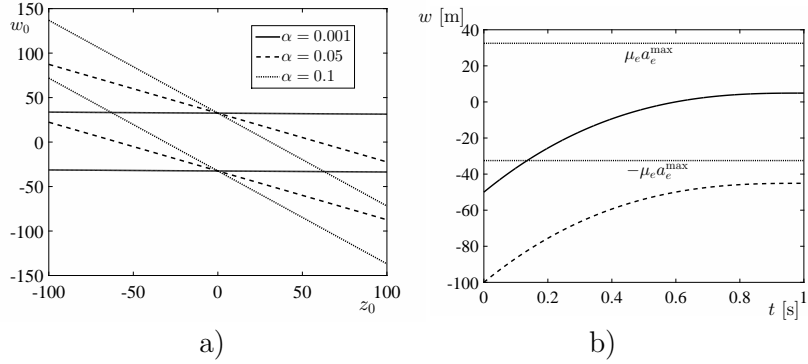

Figure 2: URG solution

In Fig. 2a, the set $\Omega$, given by (37), is shown for $\alpha = 0.001$, $\alpha = 0.05$ and $\alpha = 0.1$ ($a = 0.017$, $a = 0.92$ and $a = 1.98$, respectively). In this numerical example, $t_c = 0.9$ s, $\mu_c = 0.325$ and the terminal constraint (28) is $|w(t_f)| \leq 32.5$. In Fig. 2b, two optimal $w$-trajectories are shown for $\alpha = 0.05$ and different initial conditions. If the game starts from $(z_0 = 100, w_0 = -50) \in \Omega$, then $w(t_f) = 4.895$ m satisfies the terminal inequality constraint (the trajectory is shown by the solid line). If the initial position is $(z_0 = 100, w_0 = -100) \notin \Omega$, $w(t_f) = -45.105$ m and the terminal constraint is violated (dashed-line trajectory). The dotted lines depict the boundaries $\pm \mu_e a_e^{\max} = \pm 32.5$ m.

## 8.4 Equality Constraint Games

Consider the Equality Constraint Games, ERG$^+$ and ERG$^-$, for the same parameters as in the previous subsection, $\alpha = 0.05$, $(z_0, w_0) = (100, -100) \notin \Omega$. In this numerical example,

$$G = \begin{bmatrix} 3.72 & 2.04 \\ -2.04 & 5.91 \end{bmatrix}, \quad b^+ = \begin{bmatrix} 100 \\ -132.5 \end{bmatrix}, \quad b^- = \begin{bmatrix} 100 \\ -67.5 \end{bmatrix},$$

$$\omega_f^+ = \begin{bmatrix} z_f^+ \\ v_f^+ \end{bmatrix} = \begin{bmatrix} 32.92 \\ -11.05 \end{bmatrix}, \quad \omega_f^- = \begin{bmatrix} z_f^- \\ v_f^- \end{bmatrix} = \begin{bmatrix} 27.85 \\ -1.80 \end{bmatrix},$$

where the vectors $\omega_f^+$ and $\omega_f^-$ are defined in (75) and (86), respectively. The matrices $\tilde{G}$, $\bar{G}$ and the vectors $\chi_0$, $\gamma^+$, $\gamma^-$, defined in (76), (79), (80), (81) and (86), respectively, are

$$\tilde{G} = \begin{bmatrix} 3.72 & 2.04 \\ 2.04 & -5.91 \end{bmatrix}, \quad \bar{G} = \begin{bmatrix} 0.23 & -0.08 \\ -0.08 & -0.14 \end{bmatrix}, \quad \chi_0 = \begin{bmatrix} 100 \\ -100 \end{bmatrix}, \quad \gamma^\pm = \begin{bmatrix} 0 \\ \mp 32.5 \end{bmatrix}.$$

Thus, the game values of the ERG$^+$/ERG$^-$, given in (74), (78) and (85), respectively, are $(J^+)^* = 1821.6$, $(J^-)^* = 2659.1$.

In Fig. 3, the optimal trajectories in the games ERG$^+$ and ERG$^-$ are shown: two $z$-trajectories (Fig. 3a) and two $w$-trajectories (Fig. 3b), generated by the saddle-point pairs $(u_p^+, u_e^+)$ and $(u_p^-, u_e^-)$, respectively. It is seen that in the ERG$^+$, $w(t_f) = \mu_e a_e^{\max} = 32.5$ m, whereas in the ERG$^-$, $w(t_f) = -\mu_e a_e^{\max} = -32.5$ m, i.e., the equality terminal constraints (90) and (99) are valid, meaning that $u_e^+(\cdot) \in U_e^+$ and $u_e^-(\cdot) \in U_e^-$. Moreover, in the ERG$^+$, $z(t_f) = z_f^+ = 32.92$ m, whereas in the ERG$^-$ $z(t_f) = z_f^- = 27.85$ m, i.e. the terminal conditions (92) and (100) hold. The respective optimal controls $u_p^\pm$ and $u_e^\pm$ are shown in Figs. 4a and 4b, respectively.

In order to illustrate the saddle-point inequality (88), let us choose a constant strategy $u_e = \bar{u}_e^+ \equiv \text{const} \in U_e^+$, i.e., satisfying

$$w(t_f) = w_0 + \bar{u}_e^+ \int_0^{t_f} g_e(t) dt = \mu_e a_e^{\max}. \tag{157}$$

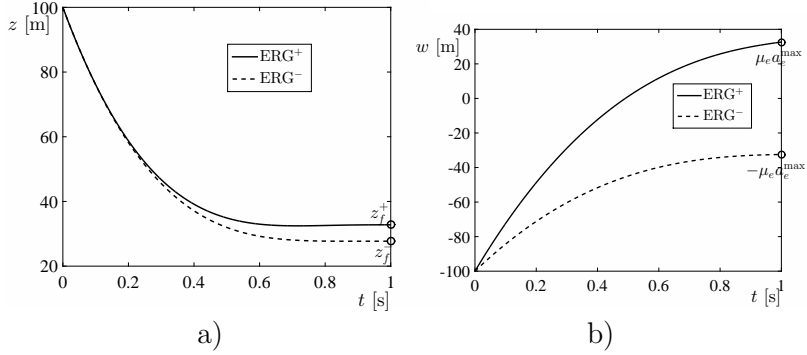

a)            b)

Figure 3: $ERG^+$ and $ERG^-$ optimal trajectories

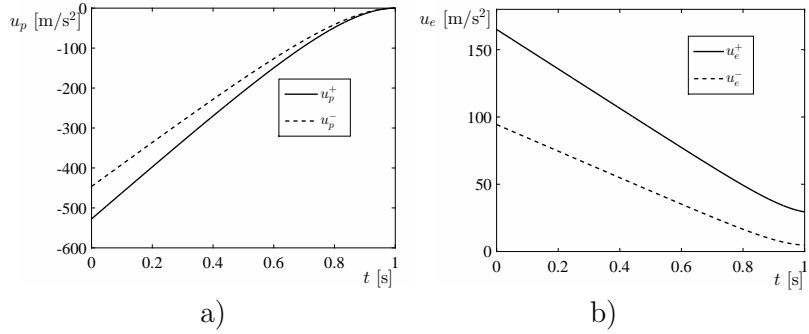

a)            b)

Figure 4: $ERG^+$ and $ERG^-$ optimal controls

For the parameters of this subsection, $\bar{u}_e^+ = 101.92$. It is calculated that $J(u_p^+(\cdot), \bar{u}_e^+) = 1358.4 < (J^+)^* = 1821.6 = J(u_p^+(\cdot), u_e^+(\cdot))$. For $u_p(t) = 400(t - t_f)$, $J(u_p(\cdot), u_e^+(\cdot)) = 2369.3 > (J^+)^*$. Thus the saddle-point inequality (88) is valid.

Due to (72), the optimal ERG$^+$-strategies $u_p^+(\cdot)$ and $u_e^+(\cdot)$ are the pointwise limits of the respective optimal UPG$_\varepsilon^+$-strategies $u_{p\varepsilon}^+(\cdot)$ and $u_{e\varepsilon}^+(\cdot)$. In Fig. 5, the difference functions $\Delta u_{p\varepsilon}^+(t) = |u_p^+(t) - u_{p\varepsilon}^+(t)|$ (Fig. 5a) and $\Delta u_{e\varepsilon}^+(t) = |u_e^+(t) - u_{e\varepsilon}^+(t)|$ (Fig. 5b) are shown for decreasing values of $\varepsilon$. The pointwise convergence is well seen.

Figs. 6 illustrate the convergence of the values $z(t_f)$ and $w(t_f)$ on the optimal trajectories in the games UPG$_\varepsilon^\pm$ to the terminal conditions $z(t_f) = z_f^\pm$ and $w(t_f) = \pm\mu_e a_e^{\max}$ in the games ERG$^\pm$. The values of $w(t_f)$ and $z(t_f)$ are shown as functions of $\varepsilon$ in Fig. 6a and Fig. 6b, respectively. In Fig.

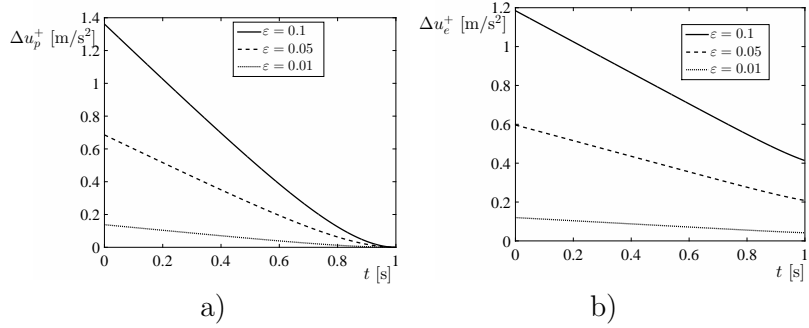

Figure 5: Convergence of $\text{UPG}_\varepsilon^+$ to $\text{ERG}^+$: optimal controls

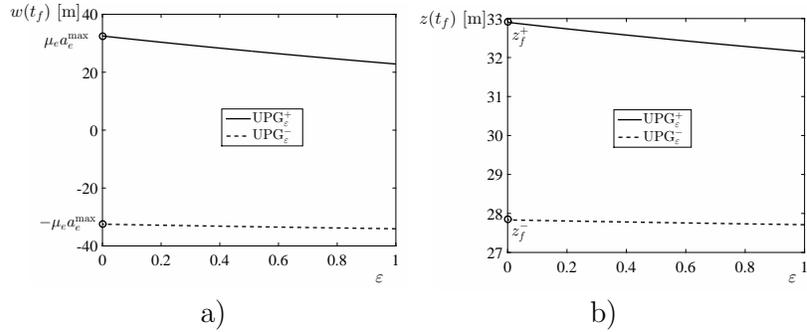

Figure 6: Convergence of $\text{UPG}_\varepsilon^\pm$ to $\text{ERG}^\pm$: terminal boundary conditions

7, the values of the games $\text{UPG}_\varepsilon^+$ and $\text{UPG}_\varepsilon^-$ are depicted as functions of $\varepsilon$, demonstrating the convergence to the values of the games $\text{ERG}^+$ and $\text{ERG}^-$, respectively.

### 8.5 Equality Constraint Game ERG

We continue using the same parameters as in the previous subsections. Remember that in this example, $a = 0.92$, $\mu_e a_e^{\max} = 32.5$. Due to Proposition 3, the pair $(u_p^+(\cdot), u_e^+(\cdot))$ constitutes the saddle point in the ERG if and only if $(z_0, w_0) \in \Omega^+ = \{w_0 + 0.92 z_0 \geq 32.5\}$ (see Fig. 8a). Similarly, for $(z_0, w_0) \in \Omega^- = \{w_0 + 0.92 z_0 \leq -32.5\}$ the saddle point in the ERG is $(u_p^-(\cdot), u_e^-(\cdot))$ (see Fig. 8b).

Let us chose the initial position $(z_0, w_0) = (100, 50) \in \Omega^+$ and calculate the cost functional (27) for three pairs of control functions: for $(u_p^+(\cdot), u_e^+(\cdot))$, for $(u_p^+(\cdot), u_e^-(\cdot))$ and for $(u_p^-(\cdot), u_e^+(\cdot))$. Similarly, for $(z_0, w_0) = (-100, -20) \in \Omega^-$ we calculate (27) for $(u_p^-(\cdot), u_e^-(\cdot))$, for $(u_p^-(\cdot), u_e^+(\cdot))$ and for $(u_p^+(\cdot), u_e^-(\cdot))$. The results are presented in Table 1.

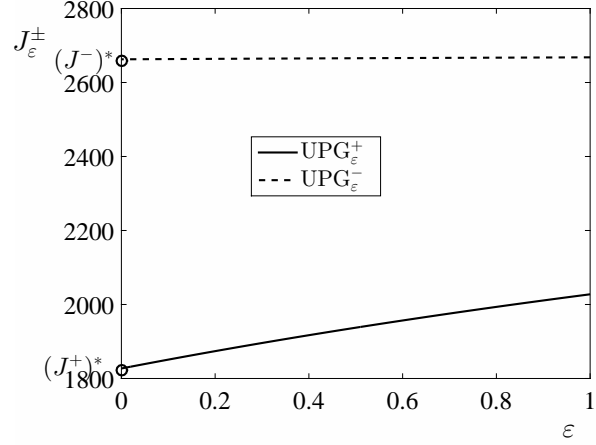

Figure 7: Convergence of $\text{UPG}_\varepsilon^+$ to $\text{ERG}^+$: game value

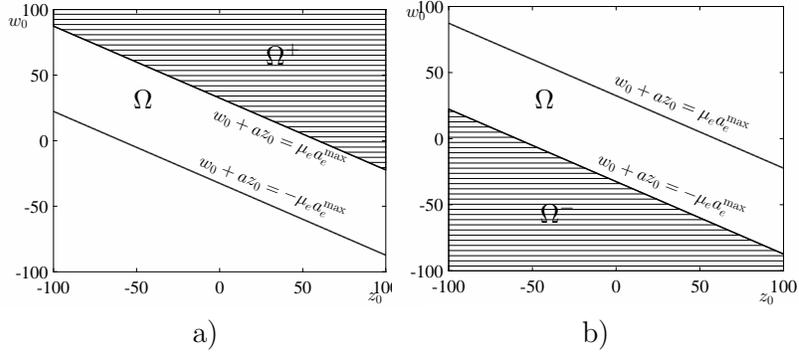

Figure 8: The sets $\Omega^+$ and $\Omega^-$

It is seen that for $(z_0, w_0) \in \Omega^+$, the saddle point inequality

$$J(u_p^+(\cdot), u_e^-(\cdot)) < J(u_p^+(\cdot), u_e^+(\cdot)) < J(u_p^-(\cdot), u_e^+(\cdot))$$

is satisfied for $(u_p^+(\cdot), u_e^+(\cdot))$. For $(z_0, w_0) \in \Omega^-$, the saddle point inequality

$$J(u_p^-(\cdot), u_e^-(\cdot)) < J(u_p^-(\cdot), u_e^-(\cdot)) < J(u_p^+(\cdot), u_e^-(\cdot))$$

Table 1: ERG Results

| $(z_0, w_0) \in \Omega^+$ | | $(z_0, w_0) \in \Omega^-$ | |
|---|---|---|---|
| Controls | Result | Controls | Result |
| $(u_p^+(\cdot), u_e^+(\cdot))$ | 1939.2 | $(u_p^-(\cdot), u_e^-(\cdot))$ | 2488.2 |
| $(u_p^+(\cdot), u_e^-(\cdot))$ | 418.8 | $(u_p^-(\cdot), u_e^+(\cdot))$ | 1463.1 |
| $(u_p^-(\cdot), u_e^+(\cdot))$ | 2347.7 | $(u_p^+(\cdot), u_e^-(\cdot))$ | 2836.7 |

is satisfied for $(u_p^-(\cdot), u_e^-(\cdot))$.

Due to Theorems 4 – 5, for $(z_0, w_0) \in \Omega^+$, the saddle point in the Reduced Game is $(u_p^+(\cdot), u_e^+(\cdot))$, whereas, for $(z_0, w_0) \in \Omega^-$, it is $(u_p^-(\cdot), u_e^-(\cdot))$.

# 9 Conclusions

A defender-attacker-target problem with non-moving target was modeled by a pursuit-evasion finite horizon linear-quadratic differential game with a terminal inequality constraint. In this game the pursuer models the defender, while the evader models the attacker. The game was solved in the following two cases: (1) the transfer functions for the players' controllers are of any orders, and (2) the transfer functions for the players' controllers are of the first-order. The cost functional of the game is a weighted sum of the square of the final range between the players, and integrals of the squares of the players' controls. The objective of the pursuer is to minimize the cost functional, while the evader has two objectives: (1) to maximize the cost functional, and (2) to hit the non-moving target. The second evader's objective yields an inequality constraint, imposed on the evader's position at the final time instant of the game.

For an unconstrained version of this game, the necessary and sufficient condition, guaranteeing the slackness of the terminal constraint, was derived. This condition yields the decomposition of the state space of the game into two non-intersecting regions, in the first of which the original game is equivalent to the unconstrained game. The open-loop saddle point solution of this game was obtained.

For the original game, considered in the second region, the condition of its equivalency to a linear-quadratic differential game with a terminal state equality constraint was obtained. The latter game is decomposed into two subgames with simpler terminal state equality constraints. For each of these subgames an open-loop saddle point solution was derived based on the penalty function method with a quadratic penalty function and a penalty coefficient tending to infinity. Using these subgames' solutions, conditions for the existence of open-loop saddle point solutions to the original game in the second region were derived and the saddle points themselves were obtained.

Being the solution in open-loop controls, these results are rather theoretic. The future issues of the topic, requiring further investigations, are (i) based on the open-loop solution, to solve the considered game in feedback strategies which is not only of a theoretical interest but also of a practical implementation; (ii) by numerical simulation, to analyse how the game parameters (especially, the penalty coefficients $\alpha$ and $\beta$) influence the control bounds.